\documentclass[12pt]{article}
\usepackage{amsmath,amssymb}

\title{Pre-Hilbert spaces without orthonormal bases}
\author{Saka\'e Fuchino\,${}^\ast$ 
}
\newif\iftesting

\newif\ifcommented 
\commentedtrue
\newif\ifJapanese
\renewcommand{\baselinestretch}{1.22}
\newcommand{\bbd}[1]{{\mathbb{#1}}}
\let\Label\label%
\iftesting
\def\label#1{\mbox{}\marginpar{{\tiny #1}}\Label{#1}\ignorespaces}%
\fi

\ifcommented
\fi
{\ifcommented\end{footnotesize}\medskip\\\fi}

\newtheorem{Thm}{{\bf Theorem}}[section]

\newtheorem{Cor}[Thm]{{\bf Corollary}}
\newtheorem{Prop}[Thm]{{\bf Proposition}}
\newtheorem{Lemma}[Thm]{{\bf Lemma}}

\newtheorem{Claim}{{\bf Claim}}[Thm]

\newtheorem{Ex}[Thm]{{\bf Example}}
\newcommand{\Thmof}[1]{{\rm Theorem \ref{#1}}}
\newcommand{\bfThmof}[1]{{\bf Theorem \ref{#1}}}

\newcommand{\Corof}[1]{{\rm Corollary \ref{#1}}}
\newcommand{\Propof}[1]{{\rm Proposition \ref{#1}}}
\newcommand{\Lemmaof}[1]{{\rm Lemma \ref{#1}}}
\newcommand{\bfLemmaof}[1]{{\bf Lemma \ref{#1}}}

\newcommand{\Claimof}[1]{{ Claim \rm\ref{#1}}}
\newcommand{\Exof}[1]{{\rm Example \ref{#1}}}

\newcommand{\Thmabove}{{\rm Theorem \number\theThm}}

\newcommand{\Claimabove}{{Claim \rm\number\theClaim}}

\newcommand{\prf}{{\bf Proof.\ }\ignorespaces}

\newcommand{\prfofClaim}{\raisebox{-.4ex}{\Large $\vdash$\ \ }}

\newcommand{\sectionof}[1]{Section {\ref{#1}}}
\newsavebox{\qedbox}\sbox{\qedbox}{
{\unitlength=0.07mm \begin{picture}(40,60)
\put(0,0){\framebox(30,44)[cc]{}}
\put(30,-7){\rule{7\unitlength}{44\unitlength}}
\put(10,-7){\rule{27\unitlength}{7\unitlength}}
\end{picture}}}
\newcommand{\qed}{\mbox{}\hfill\usebox{\qedbox}}
\newcommand{\smallqed}%
{\mbox{}\smallskip\hfill\raisebox{-.4ex}{\Large $\dashv$}\\}
\newcommand{\qedof}[1]%
{\mbox{} \hspace*{\fill}{\usebox{\qedbox}{\rm~(#1)}}%
\mbox{}}
\newcommand{\Qedof}[1]%
{\mbox{} \hspace*{\fill}{\usebox{\qedbox}%
{\rm~(#1~\number\theThm)}}}
\newcommand{\qedofThm}{\Qedof{Theorem}}

\newcommand{\qedofCor}{\Qedof{Corollary}}
\newcommand{\qedofProp}{\Qedof{Proposition}}
\newcommand{\qedofLemma}{\Qedof{Lemma}}

\newcommand{\qedofEx}{\Qedof{Example}}
\newcommand{\qedskip}{\medskip}
\newcommand{\qedofClaim}%
{\mbox{}\hfill\raisebox{-.4ex}{\Large $\dashv$ }\nolinebreak%
\mbox{\rm~({\rm Claim}~\rm\number\theClaim)}}
\newcommand{\qedofSubclaim}%
{\mbox{}\hfill\raisebox{-.4ex}{\Large $\dashv$ }\nolinebreak%
\mbox{\rm~({\rm Subclaim}~\rm\number\theSubclaim)}}

\newcommand{\st}{such that}
\newcommand{\wrt}{with respect to}

\newcommand{\wolog}{without loss of generality}
\newcommand{\Wolog}{Without loss of generality}
\newcommand{\po}{partial ordering}

\newcommand{\ctenten}{\ifmmode,\ldots\else\mbox{\rm,\,\ldots} \fi}
\newcommand{\ctentenc}{\ifmmode,\ldots,\else\mbox{\rm,\,\ldots,\,} \fi}

{\end{minipage}\end{trivlist}}
{\hfill\mbox{}\end{trivlist}}
\newcounter{frml}[section]
\def\thefrml{{\arabic{section}.\arabic{frml}}}
\def\frmlabel#1{\refstepcounter{frml}{\def\baka{#1}\ifx\baka\empty\else\Label{#1}\fi}%
{\rm({\thefrml})\hfill\hfill\hfill}}
\def\xitem[#1]{\item[\frmlabel{#1}]\mbox{}%
	\iftesting\marginpar{{\renewcommand{%
				\baselinestretch}{0.6}\tiny#1}}\fi\ignorespaces}
\def\xxitem[#1][#2]{\item[(\ref{#1}{\makebox[1.4ex][c]{#2}})]\mbox{}%
	\iftesting\marginpar{{\renewcommand{%
				\baselinestretch}{0.6}\tiny\{#1\}\{#2\}}}\fi\ignorespaces}
\def\xitemof#1{{\rm({\ref{#1}})}}

\newenvironment{xitemize}{\begin{list}{}{\parsep=0.5\smallskipamount%
			\itemindent=-0.4ex%
			\itemsep=0.5\smallskipamount\leftmargin=3.6em\labelwidth=2.73em\labelsep=0.63em}}%
							 {\end{list}}

\def\assertof#1{{\rm (#1)}}

	{\begin{trivlist}\item[\,\,#1]\mbox{}\hfill$\displaystyle}%
	{$\hfill\mbox{}\end{trivlist}}

\newcommand{\setof}[2]{\{#1\,:\,#2\}}

\newcommand{\ssetof}[1]{\{#1\}}
\newcommand{\seqof}[2]{\langle#1\,:\,#2\rangle}
\newcommand{\pairof}[1]{\langle#1\rangle}

\newcommand{\cardof}[1]{\mathopen{|\,}#1\mathclose{\,|}}
\newcommand{\mapping}[3]{#1:#2\rightarrow #3}
\newcommand{\elembed}[3]{#1:#2\stackrel{\preccurlyeq\hspace{0.8ex}}{\rightarrow}#3}
\newcommand{\id}{{\rm id}}
\newcommand{\crit}{{\rm crit}}
\newcommand{\fnsp}[2]{\mbox{}^{#1^{\mbox{}\!}}#2}

\newcommand{\osmmdin}[1]{\mathbin{\|_{#1}}}

\newcommand{\restr}{\restriction}
\newcommand{\imageof}{{}^{\,{\prime}{\prime}}}
\newcommand{\reals}{\bbd{R}}

\newcommand{\complexnrs}{\bbd{C}}

\newcommand{\poP}{\bbd{P}}

\newcommand{\regembed}%
	{\mathrel{{<}\llap{\raisebox{0.2ex}{$\scriptstyle\,\circ$}}}}
\newcommand{\forces}[2]{\,\|\hspace{-.35ex}\mbox{\sf--}_{\,#1\,}%
\mbox{\rm``}\,#2\,\mbox{\rm''}}

\newcommand{\oplusbar}{\overline{\oplus}}
\newcommand{\lessnoneq}%
{\mathrel{\raisebox{-0.8ex}{$\stackrel{<}{\scriptstyle\,\not=\,}$}}}

\newcommand{\xmbox}[1]{ $\relax{\rm #1}\relax$ }

\newcommand{\supp}{\mathop{\rm supp}}

\newcommand{\cf}{\mathop{cf\/}}

\newcommand{\otp}{\mathop{\rm otp}}
\newcommand{\cls}{\mathop{\rm cls}}
\newcommand{\Fn}{\mathop{\rm Fn}}

\newcommand{\ADS}{{\sf ADS}}
\newcommand{\FRP}{{\sf FRP}}

\newcommand{\Lim}{{\rm Lim}}

\newcommand{\calA}{{\mathcal A}}
\newcommand{\calB}{{\mathcal B}}
\newcommand{\calC}{{\mathcal C}}
\newcommand{\calD}{{\mathcal D}}

\newcommand{\calG}{{\mathcal G}}
\newcommand{\calH}{{\mathcal H}}
\newcommand{\calM}{{\mathcal M}}
\newcommand{\calN}{{\mathcal N}}
\newcommand{\calP}{{\mathcal P}}
\newcommand{\calS}{{\mathcal S}}
\newcommand{\calV}{{\mathcal V}}

\newcommand{\bba}{{\bf a}}
\newcommand{\bbb}{{\bf b}}
\newcommand{\bbc}{{\bf c}}

\newcommand{\bbe}{{\bf e}}
\newcommand{\bbf}{{\bf f}}
\newcommand{\bbg}{{\bf g}}
\newcommand{\bbu}{{\bf u}}
\newcommand{\bbv}{{\bf v}}
\newcommand{\bbx}{{\bf x}}
\newcommand{\bby}{{\bf y}}
\newcommand{\bbzero}{{\bf0}}
\newcommand{\utildeT}[1]{%
	\hbox to 0pt{$\mathop{#1}\limits_{\raise0.2ex\hbox{$\scriptstyle\sim$}}$\hss}%
		\relax\phantom{\underline{#1}}}
\newcommand{\utildeS}[1]{%
	\hbox to 0pt{$\mathop{\scriptstyle #1}\limits_{\scriptscriptstyle\sim}$\hss}%
		\relax\phantom{\underline{#1}}}
\newcommand{\utildeSS}[1]{%
	\hbox to 0pt{$\mathop{\scriptscriptstyle #1}%
		\limits_{\scriptscriptstyle\sim}$\hss}%
		\relax\phantom{\underline{#1}}}
\newcommand{\utilde}[1]{%
	\mathchoice{\utildeT{#1}}{\utildeT{#1}}{\utildeS{#1}}{\utildeSS{#1}}}
\newcount\minute	
\newcount\hour		
\newcount\hourMins  
\ifJapanese
\def\today%
{
  \the\year 年\,\zeroPadTwo{\the\month}月\,\zeroPadTwo{\the\day}日%
}
\fi
\def\now%
{
  \minute=\time    
  \hour=\time \divide \hour by 60 
  \hourMins=\hour \multiply\hourMins by 60
  \advance\minute by -\hourMins 
  \zeroPadTwo{\the\hour}:\zeroPadTwo{\the\minute}%
}
\def\zeroPadTwo#1%
{
  \ifnum #1<10 0\fi    
  #1
}
\begin{document}
\renewcommand{\thefootnote}{}
\footnotetext{{\it Date:} April 3, 2016
  \qquad {\it Last update:} 
  \today\ (\now\ JST)\vspace{-1\smallskipamount}
}
\footnotetext{\vspace{-1\smallskipamount}}
\footnotetext{\vspace{-2\smallskipamount}}   
\footnotetext{}

\maketitle
\renewcommand{\thefootnote}{$\ast$}
\footnotetext{Graduate School of System Informatics, Kobe University, 
  Kobe, Japan, \\
  \mbox{}\hspace{4.2ex}E-Mail: {\tt fuchino@diamond.kobe-u.ac.jp}\\
  The research was partially supported by Grant-in-Aid for Exploratory 
  Research No.\ 26610040 of the Ministry of Education,
  Culture, Sports, Science and Technology Japan (MEXT). 

  The author would like to thank Hiroshi Fujita of Ehime University who brought 
  the author's attention to the maximal orthonormal system in a pre-Hilbert 
  space in \Exof{example-1} which is not an orthonormal basis. Fujita also 
  suggested the 
  possibility of a construction of pre-Hilbert spaces without orthonormal bases 
  based on this example.  Thanks are also due to Hiroshi Ando of Chiba 
  University and Hiroshi Sakai of Kobe University as well as Ilijas Farah for 
  valuable comments.  
  
}
\begin{abstract}
  We give an algebraic characterization of pre-Hilbert spaces with an orthonormal 
  basis. This characterization is used to show that there are pre-Hilbert spaces
  $X$ of dimension and density $\lambda$ for any uncountable $\lambda$ without 
  any orthonormal basis. 

  Let us call a pre-Hilbert space without any orthonormal bases pathological. 
  The pair of the cardinals $\kappa\leq\lambda$ such that there is a  
  pre-Hilbert space of dimension $\kappa$ and density $\lambda$ are known to be 
  characterized  
  by the inequality $\lambda\leq\kappa^{\aleph_0}$. Our result implies that there 
  are pathological pre-Hilbert spaces with dimension $\kappa$ and 
  density $\lambda$ for all combinations of such $\kappa$ and $\lambda$ including 
  the case $\kappa=\lambda$.  

  A Singular Compactness Theorem on pathology of pre-Hilbert spaces is obtained. 

  A reflection theorem asserting that for any pathological pre-Hilbert space $X$
  there are stationarily many pathological sub-inner-product-spaces $Y$ of $X$ of 
  smaller density is shown to be equivalent with Fodor-type Reflection 
  Principle (FRP). 
\end{abstract}

\begin{quotation}
  \noindent
  \mbox{}\hfill {\bf Contents}\hfill \mbox{}\medskip\\
  \small
  \noindent
  \ref{intro}\ \ Introduction\ \ \dotfill\ \ \pageref{intro}\\
  \ref{section-2}\ \ Pathological pre-Hilbert spaces constructed from a pre-ladder 
  system\ \ \dotfill\ \ \pageref{section-2}\\ 
  \ref{char}\ \ A characterization of the non-pathology\ \ \dotfill\ \ \pageref{char}\\
  \ref{dim-den}\ \ Dimension and density of pre-Hilbert spaces\ \ \dotfill\ \ 
  \pageref{dim-den}\\  
  \ref{direct-sum}\ \ Orthogonal direct sum\ \ \dotfill\ \ \pageref{direct-sum}\\
  \ref{reflection}\ \ Reflection and non-reflection of pathology\ \ \dotfill\ \ 
  \pageref{reflection}\\ 
  \ref{singular}\ \ A Singular Compactness Theorem\ \ \dotfill\ \ 
  \pageref{singular}\\ 
  \ref{FRP}\ \ Reflection of pathology and Fodor-type Reflection 
  Principle\ \ \dotfill\ \  
  \pageref{FRP}\\ 
  \phantom{\ref{intro}}\ \   References\ \ \dotfill\ \ \pageref{literature}\\
\end{quotation}

\section{Introduction} 
\label{intro}
An inner product space whose topology is not necessarily 
complete is often called a {\it pre-Hilbert space}. 

In a pre-Hilbert space $X$, a maximal orthonormal system $S$ of $X$ does 
not necessarily span a dense subspace of $X$, that is, such $S$ does not need 
to be {\it an orthonormal basis}\/ (see \Exof{example-1} below). It is known that 
it is even possible that there is no orthonormal basis at all in some pre-Hilbert 
space (see \Lemmaof{halmos}). 
Let us call a pre-Hilbert space {\it pathological\/} if it  
does not have any orthonormal bases. If $X$ is not pathological, i.e.\ if it does 
have an 
orthonormal basis, then we say that such $X$ is {\it non-pathological}. 

By Bessel's inequality, it is easy to see that all maximal 
orthonormal system $S$ of a pre-Hilbert space $X$ has the same cardinality 
independently of whether $S$ is a basis of $X$ or not. 
This cardinality is called the {\it dimension} of the pre-Hilbert space $X$ and 
denoted by $\dim(X)$. 

In the following, we fix the scalar field $K$ of the pre-Hilbert spaces we 
consider in this paper to be $\reals$ or
$\complexnrs$ throughout.   

For an infinite set $S$, let
\begin{xitemize}
\xitem[] $\ell_2(S)=\setof{\bbu \in\fnsp{S}{K}}{\sum_{x\in S}(\bbu (x))^2
  <\infty}$,
\end{xitemize}
where $\sum_{x\in S}(\bbu (x))^2$ is defined as 
$\sup\setof{\sum_{x\in A}(\bbu (x))^2}{A\in[S]^{<\aleph_0}}$. $\ell_2(S)$ is 
endowed with a natural structure of inner product space with coordinatewise 
addition and scalar multiplication, as 
well as the inner product defined by
\begin{xitemize}
\xitem[] 
  $(\bbu,\bbv)=\sum_{x\in S}\bbu(x)\overline{\bbv(x)}$\ \ for $\bbu$, $\bbv\in\ell_2(S)$. 
\end{xitemize}
It is easy to see that $\ell_2(S)$ is a/the Hilbert space of density $\cardof{S}$. 

Note that any pre-Hilbert space $X$ of density $\lambda$ can be embedded densely into
$\ell_2(\lambda)$ as a sub-inner-product-space. Here we call a subspace $Y$ of a 
(pre-)Hilbert space $X$ a sub-inner-product-space of $X$ if $Y$ is a linear 
subspace of $X$ with the inner product which is the restriction of the inner 
product of $X$ to $Y$. 

For a pre-Hilbert space $X$ and $S\subseteq X$, we denote by $[S]_X$ the  
sub-inner-product-space of $X$ whose underlying set is the linear subspace of $X$
spanned by $S$. 

If $U$ is a subset of $\ell_2(S)$, we denote with
$\cls_{\ell_2(S)}(U)$ the topological closure of $[U]_{\ell_2(S)}$ in $\ell_2(S)$. We write 
simply $\cls(U)$ if it is clear in which $\ell_2(S)$ we are working.

For $x\in S$, let $\bbe^S_x\in\ell_2(S)$ be the standard unit vector at $x$ defined by
\begin{xitemize}
\xitem[] 
  $\bbe^S_x(y)=\delta_{x,y}$\ \ for $y\in S$. 
\end{xitemize}

For $\bba\in\ell_2(S)$, the support of $\bba$ is defined by 
\begin{xitemize}
\xitem[] 
  $\supp(\bba)=\setof{x\in S}{\bba(x)\not=0}$\ \
  ($=\setof{x\in S}{(\bba,\bbe^S_x)\not=0}$).  
\end{xitemize}
By the definition of $\ell^2(S)$, $\supp(\bba)$ is a countable subset of $S$ for 
all $\bba\in\ell^2(S)$. 

For a subset $U$ of $\ell_2(S)$ the support of $U$ is the set
$\supp(U)=\bigcup\setof{\supp(\bba)}{\bba\in U}$. 

For  $X\subseteq\ell_2(S)$ and $S'\subseteq S$, let
$X\downarrow S'=\setof{\bbu\in X}{\supp(\bbu)\subseteq S'}$. 
For $\bbu\in\ell_2(S)$, let $\bbu\downarrow S'\in\ell_2(S)$ be defined by 
\begin{xitemize}
\xitem[] 
$\left(\bbu\downarrow S'\right)(x)=\left\{\,
  \begin{array}{@{}ll}
    \bbu(x)  &\mbox{if }x\in S'\\[\jot]
    0 &\mbox{otherwise}
  \end{array}
  \right.$
\end{xitemize}
for $x\in S$. 
Note that $X\downarrow S'$ is not necessarily equal to
$\setof{\bbu\downarrow S'}{\bbu\in X}$\ \ (e.g., we have  
$X\downarrow\omega\not=\setof{\bbu\downarrow\omega}{\bbu\in X}$ where $X$ is the 
pre-Hilbert space defined in \Exof{example-1} below).

\begin{Ex}
  \label{example-1}
  Let $X$ be the sub-inner-product-space of $\ell_2(\omega+1)$ spanned by
  $\setof{\bbe^{\omega+1}_n}{n\in\omega}\cup\ssetof{\bbb}$ where
  $\bbb\in\ell_2(\omega+1)$  
  is defined by
  \begin{xitemize}
  \xitem[] $\bbb(\omega)=1$;
  \xitem[] $\bbb(n)=\frac{1}{n+2}$\ \ for $n\in\omega$.
  \end{xitemize}
  Then $\setof{\bbe^{\omega+1}_n}{n\in\omega}$ is a maximal orthonormal system in $X$ 
  but it is not a basis of $X$.  
\end{Ex}
\prf If $\setof{\bbe^{\omega+1}_n}{n\in\omega}$ were not maximal, then there 
would be an element $\bbc$ of X represented as a
linear combination of $\bbb$ and some of $\bbe^{\omega+1}_n$'s ($n\in\omega$) 
\st\ $\bbc$ is orthogonal to all $\bbe^{\omega+1}_n$, 
$n\in\omega$. However, any of such linear combinations has an infinite support 
and hence is not orthogonal to $\bbe^{\omega+1}_n$ for any $n$ in the support. 

$\setof{\bbe^{\omega+1}_n}{n\in\omega}$ is not an orthonormal basis of $X$ 
since
$\cls_{\ell_2(\omega+1)}(\setof{\bbe^{\omega+1}_n}{n\in\omega})
=\ell_2(\omega+1)\downarrow\omega\not=\ell_2(\omega+1)$. 
\qedofEx
\qedskip

For all separable pre-Hilbert spaces (including the $X$ in 
\Exof{example-1}), we can always find an orthonormal basis: suppose that $X$ is 
separable and let $\setof{\bba_n}{n\in\omega}$ be dense in $X$. Then, by 
Gram-Schmidt orthonormalization process, we can find an orthonormal system 
$\setof{\bbb_n}{n\in\omega}$ which spans the same dense sub-inner-product-space as that 
spanned by $\setof{\bba_n}{n\in\omega}$.  Thus there are no separable 
pathological pre-Hilbert spaces. 

The situation is different if we consider non-separable pre-Hilbert spaces. 

\begin{Lemma}{\rm (P.\,Halmos, see Gudder \cite{gudder})}
  \label{halmos}
  There are pre-Hilbert spaces $X$ of dimension $\aleph_0$ and 
  density $\lambda$ for any  $\aleph_0<\lambda\leq 2^{\aleph_0}$. 
\end{Lemma}

Note that a pre-Hilbert space $X$ with $\dim(X)<d(X)$ cannot have any orthonormal 
basis, that is, such a pre-Hilbert space is pathological. 

For any two pre-Hilbert spaces $X$, $Y$, the {\it orthogonal direct sum of\/ $X$ 
  and\/ $Y$} is  
the direct sum $X\oplus Y=\setof{\pairof{\bbx,\bby}}{\bbx\in X, \bby\in Y}$ of $X$ and $Y$ as 
linear spaces together with the inner product defined by 
$(\pairof{\bbx_0,\bby_0},\pairof{\bbx_1,\bby_1})=(\bbx_0,\bbx_1)+(\bby_0,\bby_1)$ 
for $\bbx_0$, $\bbx_1\in X$  
and $\bby_0$, $\bby_1\in Y$. A sub-inner-product-space $X_0$ of a pre-Hilbert 
space $X$ is an {\it orthogonal direct summand of\/ $X$} if there is a 
sub-inner-product-space $X_1$ of $X$ \st\ the mapping 
$\mapping{\varphi}{X_0\oplus X_1}{X}$; $\pairof{\bbx_0,\bbx_1}\mapsto \bbx_0+\bbx_1$ is an 
isomorphism of pre-Hilbert spaces. 
If this holds, we usually identify $X_0\oplus X_1$ with $X$ by $\varphi$ as above. 
\qedskip

\noindent
{\bf Proof of \bfLemmaof{halmos}.}  
Let $B$ be a linear basis (Hamel basis) of the linear space $\ell_2(\omega)$ extending
$\setof{\bbe^\omega_n}{n\in\omega}$. Note that $\cardof{B}=2^{\aleph_0}$ (Let 
$\calA$ be an almost disjoint family of infinite subsets of $\omega$ of 
cardinality $2^{\aleph_0}$. For each $a\in\calA$ let $\bbb_a\in\ell_2(\omega)$ be 
\st\ $\supp(\bbb_a)=a$. Then $\setof{\bbb_a}{a\in\calA}$ is a linearly 
independent subset of $\ell_2(\omega)$ of cardinality $2^{\aleph_0}$\,). 

Let $\mapping{f}{B}{
  \setof{\bbe^\lambda_\alpha}{\alpha<\lambda}\cup\ssetof{\bbzero_{\ell_2(\lambda)}}}$ be a 
surjection \st\ $f(\bbe^\omega_n)=\bbzero_{\ell_2(\lambda)}$ for all $n\in\omega$. 
Note that $f$ generates a linear mapping from the linear space $\ell_2(\omega)$ 
to a dense subspace of $\ell_2(\lambda)$. 

Let $U=\setof{\pairof{\bbb,f(\bbb)}}{\bbb\in B}$ and
$X=[U]_{\ell_2(\omega)\oplus\ell_2(\lambda)}$. Then this $X$ is as desired since 
$\setof{\pairof{\bbe^\omega_n,\bbzero}}{n\in\omega}$ is a maximal orthonormal system in 
$X$ while we have 
$\cls_{\ell_2(\omega)\oplus\ell_2(\lambda)}(X)=\ell_2(\omega)\oplus\ell_2(\lambda)$ 
and hence $d(X)=\lambda$. 
\qedof{\Lemmaof{halmos}} 
\qedskip

For sub-inner-product-spaces $X_0$, $X_1$ of a pre-Hilbert space $X$, we have
$[X_0\cup X_1]_X\cong X_0\oplus X_1$ with the isomorphism extending 
\begin{xitemize}
\xitem[] 
  $i_{X_0\cup X_1}=\setof{\pairof{\bbx_0,\pairof{\bbx_0,\bbzero}}}{\bbx\in X_0}
  \cup\setof{\pairof{\bbx_1,\pairof{\bbzero,\bbx_1}}}{\bbx_1\in X_1}$, 
\end{xitemize}
if we have 
\begin{xitemize}
\xitem[p-0] 
  $(\bbx_0,\bbx_1)=0$
  for any
  $\bbx_0\in X_0$ and $\bbx_1\in X_1$.
\end{xitemize}
Sub-inner-product-spaces $X_0$ and $X_1$ of a pre-Hilbert space $X$ with 
\xitemof{p-0} are said to be {\it orthogonal} to each other and this is denoted by
$X_0\perp X_1$.  

If $X_0$ and $X_1$ are sub-inner-product-spaces of $X$ and $X_0\perp X_1$, we 
identify $[X_0\cup X_1]_X$ with $X_0\oplus X_1$ by the isomorphism extending the
$i_{X_0\cup X_1}$ as above and write $[X_0\cup X_1]_X=X_0\oplus X_1$.

Similarly, if $X_i$, $i\in I$ are sub-inner-product-spaces of $X$ we denote 
$\oplus_{i\in I}X_i=[\bigcup_{i\in I}X_i]_X$ if $X_i$, $i\in I$ are pairwise 
orthogonal, that is, if 
we have $X_i\perp X_j$ for all distinct $i$, $j\in I$. 

For pairwise orthogonal sub-inner-product-paces $X_i$, $i\in I$ of $X$, we denote 
with $\oplusbar^X_{i\in I}X_i$ the maximal linear subspace $X'$ of $X$ \st\ $X'$ 
contains $\oplus_{i\in I}X_i$ as a dense subset of $X'$. Thus, we 
have $X=\oplusbar^X_{i\in I}X_i$ if $\oplus_{i\in I}X_i$ is dense in $X$. If it is 
clear in which $X$ we are working we drop the superscript $X$ and simply write
$\oplusbar_{i\in I}X_i$. 


An easy but very important fact for us is that
\begin{xitemize}
\xitem[orthonormal-basis] 
  if $X_i$, $i\in I$ are all 
  non-pathological with orthonormal bases $B_i$ for $X_i$, $i\in I$ and
  $X=\oplusbar_{i\in I}X_i$, then $X$ is also non-pathological with 
  the orthonormal basis $\bigcup_{i\in I}B_i$. 
\end{xitemize}

In the following we show that there are also pathological pre-Hilbert spaces $X$ with 
$\dim(X)=d(X)=\lambda$ for an uncountable $\lambda$.  For regular $\lambda$ this 
is shown in \Thmof{main-thm} and the general case in \Corof{oplus-2}. 

In \sectionof{char} we prove an algebraic characterization of pre-Hilbert spaces 
with orthonormal bases. 

In \sectionof{dim-den}, we give a proof of the theorem by Buhagiara, Chetcutib 
and Weber asserting that $\kappa\leq\lambda$ are dimension and density 
of a pre-Hilbert space if and only if $\lambda\leq\kappa^{\aleph_0}$ holds 
(see \Thmof{dim-d-0}). \Corof{oplus-2} implies that there are pathological 
pre-Hilbert spaces with $\dim(X)=\kappa$ and $d(X)=\lambda$ for all such $\kappa$ 
and $\lambda$. 

In sections \ref{reflection}, \ref{singular}, \ref{FRP} we study the 
set-theoretic reflection 
of the pathology of pre-Hilbert spaces. 

Our set-theoretic notation is quite standard. For the basic notions and notation 
in set-theory we do not explain here, the reader may consult Jech \cite{millenium-book} or 
Kunen \cite{kunen}.

\section{Pathological pre-Hilbert spaces constructed from a pre-ladder system}
\label{section-2}
For a cardinals $\lambda$, $\kappa$, let 
\begin{xitemize}
\xitem[] 
  $E^\kappa_\lambda=\setof{\alpha<\lambda}{\cf(\alpha)=\kappa}$. 
\end{xitemize}
For $E\subseteq E^\omega_\lambda$, 
$\calA=\seqof{A_\alpha}{\alpha\in E}$ is said to be a {\it ladder system} 
on $E$ if 
\begin{xitemize}
\xitem[p-1] $A_\alpha\subseteq\alpha$ for all $\alpha\in E$;
\xitem[p-1-0] 
 $A_\alpha$ is cofinal in $\alpha$ for all $\alpha\in E$; and
\xitem[p-2] $\otp(A_\alpha)=\omega$ for all $\alpha\in E$.  
\end{xitemize}
Note that, for any ladder system $\seqof{A_\alpha}{\alpha\in E}$, the sequence
$\seqof{A_\alpha}{\alpha\in E}$ is pairwise almost disjoint. We shall call a 
sequence $\seqof{A_\alpha}{\alpha\in E}$ of countable subsets 
of $\lambda$ a {\it pre-ladder system} if \xitemof{p-1} holds and \st\ it is 
pairwise almost disjoint. 

\begin{Thm} 
  \label{main-thm}
Suppose that $\kappa$ is a regular cardinal $>\omega_1$,
$E\subseteq E^\omega_\kappa$ is stationary and 
$\seqof{A_\xi}{\xi\in E}$ is a pre-ladder system \st\  
\begin{xitemize}
\xitem[m-4] 
  $A_\xi\subseteq\xi$ consists of successor ordinals for all $\xi\in E$. 
\end{xitemize}
If $\seqof{\bbu_\xi}{\xi<\kappa}$ is a sequence of elements 
of $\ell_2(\kappa)$ \st\ 
\begin{xitemize}
\xitem[m-5] $\bbu_\xi=\bbe^\kappa_\xi$ for all $\xi\in\kappa\setminus E$,
\xitem[m-6] $\supp(\bbu_\xi)=A_\xi\cup\ssetof{\xi}$ for all $\xi\in E$.
\end{xitemize}
Then, letting $U=\setof{\bbu_\xi}{\xi<\kappa}$, $X=[U]_{\ell_2(\kappa)}$ is a 
pathological pre-Hilbert space of dimension and density $\kappa$. 
\end{Thm}
\prf
We have $d(X)=\kappa$ since $\cls(X)=\ell_2(\kappa)$.
$\dim(X)\leq\dim(\ell_2(\kappa))=\kappa$ since $X$ is a sub-inner-product-space of
$\ell_2(\kappa)$ and $\dim(X)\geq\kappa$ since 
$\setof{\bbu_\alpha}{\alpha\in\kappa\setminus E}$ is an orthonormal system
$\subseteq X$ of cardinality $\kappa$. 

To show that $X$ is pathological, suppose toward a contradiction that 
$\seqof{\bbb_\xi}{\xi<\kappa}$ is an orthonormal basis of $X$. 

Let $\chi$ be a sufficiently large regular cardinal and let 
$\seqof{M_\alpha}{\alpha<\kappa}$ be a continuously increasing sequence of 
elementary submodels of $\calH(\chi)$ \st\ 
\begin{xitemize}
\xitem[m-7] $\cardof{M_\alpha}<\kappa$ for all $\alpha<\kappa$,
\xitem[m-8] $\seqof{A_\xi}{\xi\in E}$, $\seqof{\bbu_\xi}{\xi<\kappa}$,
  $\seqof{\bbb_\xi}{\xi<\kappa}\in M_0$,
\xitem[m-9] $\kappa_\alpha=\kappa\cap M_\alpha\in\kappa$ for all 
$\alpha<\kappa$ and $\seqof{\kappa_\alpha}{\alpha<\kappa}$ is a strictly 
  increasing sequence of ordinals cofinal in $\kappa$.
\end{xitemize}


For $\alpha<\kappa$, let
$H_\alpha=\ell_2(\kappa)\downarrow\kappa_\alpha$.   Note that 
$H_\alpha$ is a closed sub-inner-product-space of $\ell_2(\kappa)$ isomorphic to 
$\ell_2(\kappa_\alpha)$. 

Let $B_\alpha=\setof{\bbb_\xi}{\xi<\kappa_\alpha}$ for $\alpha<\kappa$. 
\begin{Claim}
  \label{claim-1}
  $\supp(B_\alpha)\subseteq\kappa_\alpha$ and $B_\alpha$ is an 
  orthonormal basis of $H_\alpha$. 
\end{Claim}
\prfofClaim
For $\xi<\kappa_\alpha$ $\bbb_\xi\in M_\alpha$ by \xitemof{m-8}. Hence
$\supp(\bbb_\xi)\in M_\alpha$. Since $\supp(\bbb_\xi)$ is countable it follows 
that $\supp(\bbb_\xi)\subseteq\kappa\cap M_\alpha=\kappa_\alpha$. Thus we have
$\supp(B_\alpha)\subseteq\kappa_\alpha$. 

For $\eta<\kappa_\alpha$, we have
\begin{xitemize}
\xitem[m-10] 
  $\calH(\chi)\models{}$``\,there are $A\in[\kappa]^{\aleph_0}$ and $c\in\fnsp{A}{K}$ \st\ 
  $\sum_{\xi\in A}c(\xi)\bbb_\xi=\bbu_\eta$\,''
\end{xitemize}
since $\seqof{\bbb_\xi}{\xi<\kappa}$ is an orthonormal basis. 
By \xitemof{m-8} and elementarity, it follows that
\begin{xitemize}
\xitem[m-11] 
  $M_\alpha\models{}$``\,there are $A\in[\kappa]^{\aleph_0}$ and $c\in\fnsp{A}{K}$ \st\ 
  $\sum_{\xi\in A}c(\xi)\bbb_\xi=\bbu_\eta$\,''.
\end{xitemize}

Let $A\in[\kappa]^{\aleph_0}\cap M_\alpha$ and
$c\in\fnsp{A}{\kappa}\cap M_\alpha$ be witnesses of \xitemof{m-11}. Since  
$A$ is countable we have $A\subseteq M_\alpha$. 
Thus $\bbu_\eta$ is a limit of linear combinations of elements of $B_\alpha$. 

It follows that $\cls\left([B_\alpha]_{H_\alpha}\right)
\supseteq\cls(\setof{\bbu_\xi}{\xi<\kappa_\alpha})=H_\alpha$. 
\qedofClaim\qedskip

Since $E$ is stationary, there is an
$\alpha^*<\kappa$ \st\ $\kappa_{\alpha^*}\in E$. Let $\kappa^*=\kappa_{\alpha^*}$. 
\begin{Claim}
  \label{claim-2}
  For any nonzero $\bba\in X$ represented as a linear combination of finitely 
  many elements of $U$ including (a non-zero multiple of)
  $\bbu_{\kappa*}$, there is $\xi<\kappa^*$ \st\
  $(\bba, \bbb_\xi)\not=0$. 
\end{Claim}
\prfofClaim
Suppose that 
\begin{xitemize}
\xitem[m-12] 
  $\bba=c\bbu_{\kappa*}+\sum_{\xi\in s}a_\xi\bbu_\xi+\sum_{\eta\in t}b_\eta\bbu_\eta$ 
\end{xitemize}
where $s\in[\kappa^*]^{<\aleph_0}$,
$t\in[\kappa\setminus(\kappa^*+1)]^{<\aleph_0}$ and 
$c$, $a_\xi$, $b_\eta\in K\setminus\ssetof{0}$ for $\xi\in s$ and $\eta\in t$. 
Since $\supp(\bbu_\xi)$, $\xi\in s$ are bounded subsets of $\kappa^*$ and
$\supp(\bbu_\eta)\cap\kappa^*$, $\eta\in t$ are finite, 
$\supp(\bba)\cap\kappa^*$ contains an end-segment of $A_{\kappa^*}$ and in 
particular it is non-empty. 

Thus 
$\bba\downarrow\kappa^*$ is a non-zero element of $H_{\alpha^*}$. By \Claimof{claim-1}, it 
follows that there is $\xi<\kappa^*$ \st\
$(\bba,\bbb_\xi)=(\bba\downarrow\kappa^*,\bbb_\xi)\not=0$. 
\qedofClaim\qedskip

By \Claimabove, there are no $\bba\in X$ as in the assertion of \Claimabove\ 
among $\bbb_\xi$, $\xi<\kappa$. It follows that 
$\kappa^*\not\in\bigcup\setof{\supp(\bbb_\xi)}{\xi<\kappa}$. This is a 
contradiction to the assumption that $\setof{\bbb_\xi}{\xi<\kappa}$ is an 
orthonormal basis of $X$ and hence of $\ell_2(\kappa)$. 
\qedofThm\qedskip

The construction of $X$ in \Thmabove\ can be further modified to obtain the 
following additional property of $X$: there is 
$\calS\subseteq[U]^{<\kappa}$ 
\st\ 
\begin{xitemize}
\xitem[m-2] $\calS$ is a stationary subset of $[U]^{<\kappa}$,
\xitem[m-3]  for all $A$, $B\in\calS$ with
  $A\subseteq B$, $[A]_X$ is an orthogonal direct summand of $[B]_X$. 
\end{xitemize}

For \xitemof{m-2} and \xitemof{m-2}, we can just 
start from a stationary and co-stationary $E$ and let
\begin{xitemize}
\xitem[p-14] 
  $\calS=\setof{U_\gamma}{\gamma\in\kappa\setminus E}$ 
\end{xitemize} 
where $U_\gamma=\setof{\bbu_\xi}{\xi<\gamma}$. 
Then $U$ and this $\calS$ are as desired: 
$\calS$ is a stationary subset of $[U]^{<\kappa}$ by the choice of $E$. 
For $U_{\gamma_0}$, $U_{\gamma_1}\in\calS$ 
with $\gamma_0<\gamma_1$, 
we have 
$\bbu_\xi\downarrow(\kappa\setminus\gamma_0)\in[U_{\gamma_1}]_X$ for all
$\xi\in\gamma_1\setminus\gamma_0$.  
Hence 
\begin{xitemize}
\xitem[] 
  $[U_{\gamma_1}]_X
  =[U_{\gamma_0}]_X
  \oplus[\setof{\bbu_\xi\downarrow(\kappa\setminus\gamma_0)}{
      \xi\in\gamma_1\setminus\gamma_0}]_X$. 
\end{xitemize}


\Thmof{main-thm} applied to $\kappa=\omega_1$ gives pathological pre-Hilbert 
spaces with interesting properties.  Note that for a stationary subset $E$ of 
$\omega_1$ there is a \po\ which ``shoots'' a club subset inside $E$ while 
preserving all cardinals (e.g. the shooting a club forcing with finite 
conditions). 

If $X$ is a pre-Hilbert space constructed as in \Thmof{main-thm} for stationary 
and co-stationary $E\subseteq E^\omega_{\omega_1}$ and a pre-ladder system 
on $E$, letting $U\subseteq\ell_2(\omega_1)$ be the generator of $X$ as in 
\Thmof{main-thm}, we have that 
$X\downarrow \alpha$ is non-pathological for all $\alpha<\omega_1$ since 
$X\downarrow\alpha$ is separable. 
If we shoot a club subset of $\omega_1\setminus E$, 
we obtain a continuously increasing sequence of non-pathological 
sub-inner-product-spaces $\seqof{X_\alpha}{\alpha<\omega_1}$ of $X$ \st\ 
$\bigcup_{\alpha<\omega_1} X_\alpha=X$ and that $X_\alpha$ is an orthogonal 
direct summand of $X_{\alpha+1}$ for all $\alpha<\omega_1$. It follows that $X$ 
is non pathological in such a generic extension. Thus we obtain: 
\begin{Cor}
  \label{cor-1}
  \assertof{1} There is a pathological pre-Hilbert space $X$ of dimension and 
  density $\aleph_1$ \st\ there is a \po\ $\poP$ preserving all cardinals \st\
  $\forces{\poP}{X\xmbox{ has an ortho\-normal basis}}$.\smallskip

  \assertof{2} There is a pathological pre-Hilbert space $X$ of dimension and 
  density $\aleph_1$ \st, for any \po\ $\poP$ preserving $\omega_1$, we 
  have $\forces{\poP}{X\mbox{ is pathological\/}}$.
\end{Cor}
\prf A proof of \assertof{1} is already explained above. 
For \assertof{2}, we can use the club set $E=E^\omega_{\omega_1}$ in the 
construction of the  
proof of \Thmof{main-thm}. The pre-Hilbert space $X$ constructed in the proof of 
\Thmof{main-thm} with this $E$ is as desired: since $E^*$ remains stationary in 
any generic extension preserving $\omega_1$, $X$ remains pathological there. 
\qedofCor 
\qedskip

\section{A Characterization of the non-pathology}
\label{char}
Using some of the ideas in the proof of \Thmof{main-thm}, we obtain an 
``algebraic'' characterization of pre-Hilbert spaces with 
orthonormal bases  (see \Thmof{T-char}). 
This characterization is used in later sections.  
\begin{Lemma}
  \label{L-char-1}
  Suppose that $X$ is a pre-Hilbert space and $X$ is a dense sub-inner-product-space 
  of $\ell_2(S)$. If $\calB\subseteq X$ is an orthonormal basis then, for any
  $S_0\subseteq S$, there is an $A\subseteq S$ \st\
  $S_0\subseteq A$, $\cardof{A}=\cardof{S_0}+\aleph_0$, 
  $X\downarrow A$ is a dense sub-inner-product-space of $\ell_2(S)\downarrow A$, 
  $\calB_A=\setof{\bbb\in\calB}{\supp(\bbb)\subseteq A}$ is an orthonormal basis 
  of $X\downarrow A$ and $\calB^-_A=\calB\setminus\calB_A$ is an orthonormal 
  basis of $X\downarrow (S\setminus A)$. In particular, we have
  $X=(X\downarrow A)\oplus(X\downarrow (S\setminus A))$. 
\end{Lemma}
\prf Let $\chi$ be a sufficiently large regular cardinal and let
$M\prec\calH(\chi)$ be \st\ 
\begin{xitemize}
\xitem[chr-1] 
  $K$, $X$, $S$, $\calB\in M$, $S_0\subseteq M$ and $\cardof{M}=\cardof{S_0}+\aleph_0$.
\end{xitemize}

We show that $A=S\cap M$ is as desired. Since $S_0\subseteq M$, we have 
$S_0\subseteq S\cap M=A$. 

Since $\calB$ is also an orthonormal basis of $\ell_2(S)$,
we have 
\begin{xitemize}
\xitem[chr-1-0] 
  $\calH(\chi)\models$``there is a $B\in[\calB]^{\aleph_0}$ and $c\in\fnsp{B}{K}$ 
  \st\ $\sum_{\bbu\in B}c(\bbu)\bbu=\bbe^S_s$\,''
\end{xitemize}
for all $s\in A$.  
By elementarity, it follows that 
\begin{xitemize}
\xitem[chr-2] 
  $M\models$``there is a $B\in[\calB]^{\aleph_0}$ and $c\in\fnsp{B}{K}$ 
  \st\ $\sum_{\bbu\in B}c(\bbu)\bbu=\bbe^S_s$\,''.
\end{xitemize}
Let $B\in[\calB]^{\aleph_0}\cap M$ and $c\in\fnsp{B}{K}\cap M$ be witnesses 
of \xitemof{chr-2}. By $B\in M$ and since $B$ is countable, we have $B\subseteq M$. 
For each $\bbb\in B$, since $\bbb\in M$ and $\supp(\bbb)$ is countable, we have
$\supp(\bbb)\subseteq M$. 
It follows that $B\subseteq\calB_A$ and 
$\bbe^S_s\in\cls_{\ell_2(S)\downarrow A}(\calB_A)$ for all $s\in A$. 

Thus 
$\setof{\bbe^S_s}{s\in A}\subseteq
\cls_{\ell_2(S)\downarrow A}(\calB_A)$ 
and hence 
\begin{xitemize}
\xitem[chr-3] 
  $\cls_{\ell_2(S)\downarrow A}(\calB_A)=\ell_2(S)\downarrow A$. 
\end{xitemize}
Since $\calB_A\subseteq X\downarrow A$, \xitemof{chr-3} implies that
$X\downarrow A$ is dense in $\ell_2(S)\downarrow A$.

For any $\bbb\in\calB^-_A$, we have $\supp(\bbb)\subseteq S\setminus A$: 
otherwise, $\bbb\downarrow A\not=\bbzero_{\ell_2(S)}$. By \xitemof{chr-3} it 
follows that there is a $\bbc\in\calB_A$ \st\
$(\bbb,\bbc)=(\bbb\downarrow A,\bbc)\not=0$. This is a contradiction to the 
orthonormality of $\calB$. 
Thus $\calB^-_A\subseteq X\downarrow(S\setminus A)$. 

$\calB^-_A$ is an 
orthonormal basis of $X\downarrow(S\setminus A)$: similarly to the argument 
above, it is enough to show  
that,  for each $s\in S\setminus A$, $\bbe^S_s$ can be obtained as a (possibly 
infinite) sum of elements in $\calB^-_A$ in $\ell_2(S)\downarrow(S\setminus A)$. 
Since $\calB$ is an orthonormal basis  
of $\ell_2(S)$, 
we have $\bbe^S_s=\sum_{\bbb\in B}(\bbe^S_s,\bbb)\bbb$ where
$B=\setof{\bbb\in\calB}{(\bbe^S_s,\bbb)\not=0}$. Since
$\supp(\bbb)\not\subseteq A$ for all $\bbb\in B$, we have $B\subseteq\calB^-_A$. 
\qedofLemma

\begin{Lemma}
  \label{L-char-2}
  Suppose that $X$ is a non-pathological pre-Hilbert space and $X$ is a dense 
  sub-inner-product space of $\ell_2(S)$ for some infinite set $S$. Then there is a 
  partition $\calP$ of  
  $S$ into countable subsets \st\ $X=\oplusbar_{A\in\calP}X\downarrow A$. 
\end{Lemma}
\prf Let $\cardof{S}=\kappa$ and $\calB=\setof{\bbb_\alpha}{\alpha<\kappa}$ be an 
orthonormal basis  
of $X$. Let $S=\setof{s_\alpha}{\alpha<\kappa}$. 

We define by induction on $\alpha\in\kappa$ the sequences 
$\seqof{S_\alpha}{\alpha<\kappa}$ and $\seqof{A_\alpha}{\alpha<\kappa}$ of 
subsets of $S$ \st:
\begin{xitemize}
\xitem[chr-4] $S_0=S$;
\xitem[chr-5] $A_\alpha\in[S_\alpha]^{\aleph_0}$ for all $\alpha\in\kappa$;
\xitem[chr-6] $S_{\alpha+1}=S_\alpha\setminus A_\alpha$ for all $\alpha\in\kappa$; 
\xitem[chr-7] $S_\gamma=\bigcap_{\alpha<\gamma}S_\alpha$ for all limit
  $\gamma\in\kappa$;
\xitem[chr-8] $s_\alpha\in\bigcup_{\beta\leq\alpha}A_\alpha$ for 
  all $\alpha\in\kappa$; 
\xitem[chr-9] $\calB\cap(X\downarrow S_\alpha)$ is an orthonormal basis of
  $X\downarrow S_\alpha$ for 
  all $\alpha\in\kappa$; and 
\xitem[chr-10] $\calB\cap(X\downarrow A_\alpha)$ is an orthonormal basis of
  $X\downarrow A_\alpha$ for 
  all $\alpha\in\kappa$. 
\end{xitemize}

The construction of $A_\alpha$ and $S_{\alpha+1}$ is possible by \Lemmaof{L-char-1}. 
We just have to check that the construction of $S_\gamma$ at limit steps 
$\gamma<\kappa$ works. 

For a limit $\gamma<\kappa$ we have $S_\gamma=\bigcap_{\alpha<\gamma}S_\alpha$ by 
\xitemof{chr-7}. For each $s\in S_\gamma$ and $\alpha<\gamma$ there are a countable
$B_\alpha\subseteq\calB\cap(X\downarrow S_\alpha)$ and a sequence
$\seqof{a^\alpha_\bbb}{\bbb\in B_\alpha}$ in $K$ \st\ 
\begin{xitemize}
\xitem[chr-11] 
  $\bbe^S_s=\sum_{\bbb\in B_\alpha}a^\alpha_\bbb\bbb$. 
\end{xitemize}
By the uniqueness of the representation of elements of $\ell_2(S)$ as an infinite 
linear combination of elements of $\calB$. It follows that there is a countable 
$B^*\subseteq\calB\cap (X\downarrow S_\alpha)$ and a sequence
$\seqof{a_\bbb}{\bbb\in B^*}$ \st\  
$B_\alpha=B^*$ for all $\alpha<\gamma$ and $a^\alpha_\bbb=a_\bbb$ for all 
$\alpha<\gamma$ and $\bbb\in B^*$. It follows that
$B^*\subseteq\calB\cap(X\downarrow B_\gamma)$.

Thus, we have 
$\bbe^S_s\in\cls_{\ell_2(S)}[\calB\cap(X\downarrow S_\gamma)]$, for all 
$s\in S_\gamma$. It follows that  
$\calB\cap(X\downarrow S_\gamma)$ is an orthonormal basis of
$X\downarrow S_\gamma$, i.e.\ $S_\gamma$ satisfies \xitemof{chr-9}. 

$\calP=\setof{A_\alpha}{\alpha<\kappa}$ is then a partition of $S$ as desired. \qedofLemma
\begin{Thm}
  \label{T-char}
  Suppose that $X$ is a pre-Hilbert space. Then $X$ is non-pathological if and only 
  if there are separable sub-inner-product-spaces $X_\alpha$, $\alpha<\delta$ of
  $X$ \st\ $X=\oplusbar_{\alpha<\delta}X_\alpha$.  
\end{Thm}
\prf If $X$ is separable then the claim is trivial with $\delta=1$. 

Suppose that $X$ is non-separable.

If $X$ is non-pathological then there are separable 
sub-inner-product-spaces $X_\alpha$, $\alpha<\kappa$ for $\kappa=d(X)$ with
$X=\oplusbar_{\alpha<\kappa}X_\alpha$ by \Lemmaof{L-char-2}. 

Conversely, if there are $X_\alpha$, $\alpha<\delta$ as above, then 
each $X_\alpha$ for $\alpha\in\delta$ has an orthonormal basis $B_\alpha$. 
$\calB=\bigcup_{\alpha<\delta}B_\alpha$ is then an orthonormal basis of $X$. \qedofThm

\begin{Lemma}
  \label{filtration}
  Suppose that $X$ is a non-pathological pre-Hilbert space and $X$ is a dense 
  sub-inner-product space of $\ell_2(S)$ for some uncountable set $S$. Then there 
  is a filtration $\seqof{S_\alpha}{\alpha<\kappa}$ of $S$ for
  $\kappa=\cf(\cardof{S})$ \st\ 
  $X\downarrow S_\alpha$ is an orthogonal direct summand of $X$ for all
  $\alpha<\kappa$. 
\end{Lemma}
\prf By \Lemmaof{L-char-2} there is a partition $\calP$ of $S$ into countable 
subsets \st\ $X=\oplusbar_{P\in\calP}X\downarrow P$. Let 
$\seqof{\calP_\alpha}{\alpha<\kappa}$ be a filtration of $\calP$ and let 
$S_\alpha=\bigcup\calP_\alpha$ for $\alpha<\kappa$. Then 
$\seqof{S_\alpha}{\alpha<\kappa}$ is as desired. \qedofLemma
\qedskip

The following Lemmas are used in \sectionof{FRP}. We put them together here since they 
stand in a similar context as that of previous results in this section.

\begin{Lemma}
  \label{summand}
  Suppose that $X$ is a pre-Hilbert-space which is a dense 
  sub-inner-product-space of $\ell_2(S)$. For $S'\subseteq S$ \st\ 
  \begin{xitemize}
  \xitem[chr-15-0] 
    $X\downarrow S'$ is dense in $\ell_2\downarrow S'$,
  \end{xitemize}
  $X\downarrow S'$ is not an 
  orthogonal direct summand of $X$ if and only if there is $\bba\in X$ \st\
  \begin{xitemize}
  \xitem[chr-16] 
    $\bba\downarrow S'\not\in X$. 
  \end{xitemize}
\end{Lemma}
\prf If there is no $\bba\in X$ with \xitemof{chr-16} then we clearly have
$X=(X\downarrow S')\otimes(X\downarrow S\setminus S')$. 

Suppose that $\bba\in X$ satisfies \xitemof{chr-16}. Note that then we have
$\supp(\bba)\nsubseteq S'$ and $\supp(\bba)\cap S'\not=\emptyset$. 
Suppose toward a contradiction that there is a sub-inner-product space $X''$ of 
$X$ \st\ 
\begin{xitemize}
\xitem[chr-17] 
  $X=(X\downarrow S')\oplus X''$. 
\end{xitemize}
Then there are $\bba'\in X\downarrow S'$ 
and $\bba''\in X''$ \st\ $\bba=\bba'+\bba''$. So $\bba''=\bba-\bba'$. It follows 
that $\bba''\downarrow S'\not=\bbzero$ by \xitemof{chr-16}. By \xitemof{chr-15-0}, 
there is some $\bbb\in X\downarrow S'$ \st\
$(\bba'',\bbb)=(\bba''\downarrow S',\bbb)\not=0$. This is a contradiction to 
\xitemof{chr-17}.\qedofLemma

\begin{Lemma}
  \label{L-dense-0}
  Suppose that $X$ is a pre-Hilbert-space which is a dense 
  sub-inner-product-space of $\ell_2(S)$. For a sufficiently large regular $\chi$ 
  and $M\prec\calH(\chi)$ with $K$, $X$, $S\in M$ , $X\downarrow (S\cap M)$ is 
  dense in $\ell_2(S)\downarrow(S\cap M)$. 
\end{Lemma}
\prf For $s\in S\cap M$, we have
\begin{xitemize}
\xitem[chr-18] $\calH(\chi)\models$ there are $A\in[X]^{\aleph_0}$ and $c\in\fnsp{A}{K}$
  \st\ $\bbe^S_s=\sum_{\bbb\in A}c(\bbb)\bbb$. 
\end{xitemize}
By elementarity it follows that 
\begin{xitemize}
\xitem[chr-19] $M\models$ there are $A\in[X]^{\aleph_0}$ and $c\in\fnsp{A}{K}$
  \st\ $\bbe^S_s=\sum_{\bbb\in A}c(\bbb)\bbb$. 
\end{xitemize}
Let $A\in[X]^{\aleph_0}\cap M$ and $c\in\fnsp{A}{K}\cap M$ be witnesses of 
\xitemof{chr-19}. By the countability of $A$ we have $A\subseteq M$ and, for each
$\bbb\in A$, $\supp(\bbb)\subseteq M$ since $\supp(\bbb)$ is countable. 

This shows that $\bbe^S_s\in\cls(X\downarrow (S\cap M))$. 
\qedofLemma
\begin{Lemma}
  \label{L-dense-1}
  Suppose that $X$ is a pre-Hilbert-space which is a dense 
  sub-inner-product-space of $\ell_2(S)$ for an uncountable S. Then there is a 
  filtration $\seqof{S_\alpha}{\alpha<\kappa}$ of $S$ \st\
  $X\downarrow S_\alpha$ dense in $\ell_2(S)\downarrow S_\alpha$ for all $\alpha<\kappa$
\end{Lemma}
\prf Let $\chi$ be a sufficiently large regular cardinal. Let
$\kappa=\cf(\cardof{S})$ and let $\seqof{M_\alpha}{\alpha<\kappa}$ be a 
continuously increasing sequence of elementary submodels of $\calH(\chi)$ \st\ 
\begin{xitemize}
\xitem[chr-20] $K$, $X$, $S\in M_0$,
\xitem[chr-21] $\cardof{M_\alpha}<\cardof{S}$ for all $\alpha<\kappa$,
\xitem[chr-22] $S\subseteq \bigcup_{\alpha<\kappa}M_\alpha$. 
\end{xitemize}
Letting $S_\alpha=S\cap M_\alpha$ for $\alpha<\kappa$, the sequence 
$\seqof{S_\alpha}{\alpha<\kappa}$ is as desired by \Lemmaof{L-dense-0}. 
\qedofLemma

\section{Dimension and density of pre-Hilbert spaces}
\label{dim-den}
The proof of \Lemmaof{halmos} actually yields pre-Hilbert spaces 
of the following combinations of dimension and density:

\begin{Lemma}[A generalization of \bfLemmaof{halmos}]
  \label{halmos-0}
  For any cardinal $\kappa$ and $\lambda$ with 
  $\kappa<\lambda\leq\kappa^{\aleph_0}$, there are (pathological) pre-Hilbert spaces of 
  dimension $\kappa$ and density $\lambda$.\qed
\end{Lemma}

On the other hand if $\kappa^{\aleph_0}<\lambda$ there are no pre-Hilbert space 
$X$ with dimension $\kappa$ and density $\lambda$. 
\begin{Prop}{\rm (David Buhagiara, Emmanuel Chetcutib and Hans Weber 
    \cite{buhagiara-etal}, see also \cite{farah})}
  \label{dim-d}
  For any pre-Hilbert space $X$, we have
  $d(X)\leq\cardof{X}\leq (\dim(X))^{\aleph_0}$. 
\end{Prop}
\prf Let $X$ be a pre-Hilbert space. We may assume \wolog\ that $X$ is a dense 
sub-inner-product-space of the Hilbert space $\ell_2(\kappa)$ for
$\kappa=d(X)>\dim(X)\geq\aleph_0$.   

Let $\calB=\seqof{\bbb_\xi}{\xi<\kappa}$ be a maximal orthonormal system in $X$
and
$D=\bigcup\setof{\supp(\bbb_\xi)}{\xi<\kappa}$. By the assumption we have
$\cardof{D}=\kappa$. 
\begin{Claim}
  For any distinct $\bba_0$, $\bba_1\in X$ we have
  $\bba_0\restr D\not=\bba_1\restr D$.
\end{Claim}
\prf Suppose that there were $\bba_0$, $\bba_1\in X$ \st\ $\bba_0\not=\bba_1$ but
$\bba_0\restr D=\bba_1\restr D$. Then $\bba_2=\bba_1-\bba_0$ would be a non-zero element 
of $X$ orthogonal to all $\bbb_\xi$, $\xi<\kappa$. This is a contradiction to the 
maximality of $\calB$.\qedofClaim\qedskip

Let $\mapping{\varphi}{\ell_2(D)}{X}$ be defined by
\begin{xitemize}
\xitem[m-13] $\varphi(\bbc)=\left\{\,
  \begin{array}{@{}ll}
    \mbox{the unique }\bba\in X\mbox{ \st\ }\bbc=\bba\restr D; 
    &\mbox{if there is such }\bba\in X,\\[\jot]
    \bbzero; &\mbox{otherwise}
  \end{array}
\right.$
\end{xitemize}
for $\bbc\in\ell_2(D)$. $\varphi$ is well-defined by \Claimabove\ and it is 
surjective. Thus we have
\begin{xitemize}
\xitem[] 
  $d(X)\leq\cardof{X}\leq \cardof{\ell_2(D)}=(\dim(X))^{\aleph_0}$.\\\qedofProp
\end{xitemize}

The following theorem will be yet extended in \Corof{cor-3}. 

\begin{Thm}
  \label{dim-d-0}
  For any cardinal $\kappa\leq\lambda$ there is a pre-Hilbert space of dimension 
$\kappa$ and density $\lambda$ if and only if $\lambda\leq \kappa^{\aleph_0}$ holds.
\end{Thm}
\prf For  $\kappa=\lambda$, $\ell_2(\kappa)$ is an example of pre-Hilbert space 
of dimension and density $\kappa$ and $\lambda$. If
$\kappa<\lambda<\kappa^{\aleph_0}$, \Lemmaof{halmos-0} provides an example.  

The converse also holds by \Propof{dim-d}.\qedofThm

\section{Orthogonal direct sum} 
\label{direct-sum}
In a variety $\calV$ of algebraic structures it can happen that there is a non 
free algebra $A\in\calV$ \st\ the product $A\otimes F$ is free for some free algebra
$F\in\calV$. For example, it is known that there are non-free projective algebra $B$
in the variety $\calB$ of Boolean algebras but free product $B\oplus F$ of any 
projective algebra $B$ with a sufficiently large free Boolean algebra $F$ is 
free. 

In contrast, the pathology of pre-Hilbert space remains by orthogonal direct sum.
\begin{Thm}
  \label{T-direct-sum}
  For any pre-Hilbert spaces $X_0$ and $X_1$, the orthogonal direct sum
  $X_0\oplus X_1$ is pathological if and only if at least one of $X_0$ and $X_1$ 
  is pathological. 
\end{Thm}
\prf If $X_0$ and $X_1$ are both non-pathological and $\calB_0$ and $\calB_1$ are 
orthonormal bases of $X_0$ and $X_1$ respectively, then 
$\calB_0\times\ssetof{\bbzero_{X_2}}\cup\ssetof{\bbzero_{X_1}}\times\calB_1$ is an 
orthonormal basis of $X_0\oplus X_1$. 

Conversely, suppose that $X_0\oplus X_1$ is non-pathological and $\calB$ is an 
orthonormal basis of $X=X_0\oplus X_1$. \Wolog, we may assume that there are $S$,
$S^0$, $S^1$ \st\ $S=S^0\cup S^1$, $S^0\cap S^1=\emptyset$, $X_i$ is a dense 
sub-inner-product-space of $\ell_2(S)\downarrow S^i$ for $i\in 2$ and
$X_0\oplus X_1=[X_0\cup X_1]_{\ell_2(S)}$. 

By \Lemmaof{L-char-2}, there is a partition $\seqof{A_\alpha}{\alpha<\delta}$ of 
$S$ into countable sets \st\ $X=\oplusbar_{\alpha\in\kappa}X\downarrow A_\alpha$. We 
may assume that the elements of partition $A_\alpha$ in the proof of 
\Lemmaof{L-char-2} is obtained in the construction as the intersection of 
$S_\alpha$ (in the proof of \Lemmaof{L-char-2}) and countable
$M_\alpha\prec\calH(\chi)$ \st\ $\calB$, $X_0$, $X_1$, $S^0$, 
$S^1\ctenten\in M_\alpha$. 
Then as in the proof of \Lemmaof{L-char-1}, we have
$X\downarrow A_\alpha= (X_0\downarrow(A_{0,\alpha}))\oplus (X_1\downarrow(A_{1,\alpha}))$   
where $A_{i,\alpha}=A_\alpha\cap S^i$ for $i\in 2$. 

Let
$\calP_i=\setof{A_{i,\alpha}}{\alpha<\kappa,\,A_{i,\alpha}\not=\emptyset}$
for $i\in 2$. Then $X_i=\oplusbar_{P\in\calP_i}X_i\downarrow P$ for
$i\in 2$. Thus $X_i$, $i\in 2$ are non-pathological. 
\qedofThm

\begin{Cor}
  \label{oplus-2}
  For any uncountable cardinal $\lambda$, there is a pathological pre-Hilbert space $Z$ 
  of dimension and density $\lambda$.
\end{Cor}
\prf Let $X$ be any pathological pre-Hilbert space with density $\aleph_1$.  
Then $Z=X\oplus\ell_2(\lambda)$ has 
dimension and density $\lambda$. $Z$ is pathological by \Thmof{T-direct-sum}. \qedofCor

\begin{Cor}
  \label{cor-2}
  For any infinite cardinals $\kappa$ and $\lambda$ with
  $\kappa\leq \lambda\leq \kappa^{\aleph_0}$ there is a pathological pre-Hilbert space of 
  dimension $\kappa$ and density $\lambda$.
\end{Cor}
\prf By \Lemmaof{halmos-0} and \Corof{oplus-2}.\qedofCor

\begin{Cor}
  \label{cor-3}
  \assertof{1}
  For any infinite cardinals $\kappa$ and $\lambda$ with
  $\kappa\leq \lambda\leq \kappa^{\aleph_0}$ there is a pathological pre-Hilbert 
  space of  
  dimension $\kappa$ and density $\lambda$ \st\ there is a \po\ $\poP$ preserving 
  all cardinals \st\ $\forces{\poP}{X\mbox{ is non-pathological\/}}$. \smallskip

  \assertof{2}   For any infinite cardinals $\kappa$ and $\lambda$ with
  $\kappa\leq \lambda\leq \kappa^{\aleph_0}$ there is a pathological pre-Hilbert 
  space of  
  dimension $\kappa$ and density $\lambda$ which remains pathological in any 
  generic extension preserving $\omega_1$.
\end{Cor}
\prf The pre-Hilbert space of the form $X\oplus Y$ will do where $X$ is as in 
\Corof{cor-1},\assertof{1} or \assertof{2} and $Y$ is as in 
\Corof{cor-2}.\qedofCor

\section{Reflection and non-reflection of pathology}
\label{reflection}
For any pre-Hilbert space $X$ all sub-inner-product-spaces of $X$ of density
$\aleph_0$ are non-pathological. If $S\subseteq E^\omega_{\omega_2}$ is 
non-reflecting stationary set, then the sub-inner-product-space of
$\ell_2(\omega_2)$ constructed from a ladder system on $S$,  
there are club many $\beta<\omega_2$ \st\ 
$X\downarrow \beta$ is non-pathological. 

A similar non-reflection theorem holds at an arbitrary regular uncountable 
cardinal $\kappa>\aleph_1$ under a weak form of the square principle at $\kappa$.

For a regular cardinal $\kappa$, $\ADS^-(\kappa)$ is the assertion
that there is a stationary set $S\subseteq E^\omega_\kappa$ and a sequence
$\seqof{A_\alpha}{\alpha\in S}$ \st
\begin{xitemize}
\xitem[ads-0] $A_\alpha\subseteq\alpha$ and $\otp(A_\alpha)=\omega$ for all
  $\alpha\in S$;
\xitem[ads-1] for any $\beta<\kappa$, there is a mapping 
$\mapping{f}{S\cap\beta}{\beta}$ \st\ $f(\alpha)<\sup(A_\alpha)$ for all
  $\alpha\in S\cap\beta$ and $A_\alpha\setminus f(\alpha)$,
  $\alpha\in S\cap\beta$ are pairwise disjoint 
\end{xitemize}
(for more about $\ADS^-(\kappa)$, see Fuchino, Juha\'asz, Soukup, 
Szentmikl\'ossy, Usuba \cite{fuchino-juhasz-etal} and Fuchino, Sakai, Soukup \cite{more} ). 

We shall call $\seqof{A_\alpha}{\alpha\in S}$ as above 
an $\ADS^-(\kappa)$-sequence. Note that it follows from \xitemof{ads-0} and 
\xitemof{ads-1} that $A_\alpha$, $\alpha\in S$ are pairwise almost disjoint. 

Under $\ADS^-(\kappa)$, we may further assume that the $\ADS^-(\kappa)$-sequence 
$\seqof{A_\alpha}{\alpha\in S}$ satisfies that
$A_\alpha\subseteq\alpha\setminus\Lim$ for all $\alpha\in S$. 

Since an $\ADS^-(\kappa)$-sequence is a pre-ladder system, we can apply the 
construction of pre-Hilbert spaces in the proof of \Thmof{main-thm} to the 
sequence and obtain the following:

\begin{Thm}
  \label{T-refl-1}
  Assume that $\ADS^-(\kappa)$ holds for a regular cardinal $\kappa>\omega_1$. 
  Then there is a pathological dense sub-inner-product-space $X$ of $\ell_2(\kappa)$ 
  \st\ $X\downarrow \beta$ is non-pathological for all $\beta<\kappa$. Furthermore 
  for any regular $\lambda<\kappa$,
  $\setof{S\in[\kappa]^\lambda}{X\downarrow S\mbox{ is non-pathological\/}}$ 
  contains a club subset of $[\kappa]^\lambda$.
\end{Thm}
\prf
Let $\seqof{A_\alpha}{\alpha\in E}$ be an $\ADS^-(\kappa)$-sequence on a 
stationary $E\subseteq E^\omega_\kappa$. Let 
$\seqof{\bbu_\xi}{\xi<\kappa}$ be a sequence of elements of $\ell_2(\kappa)$ with 
\xitemof{m-5} and \xitemof{m-6}, $U=\setof{\bbu_\xi}{\xi<\kappa}$ and 
$X=[U]_{\ell_2(\kappa)}$. Then $X$ is pathological by \Thmof{main-thm}. 

For $\beta<\kappa$ let $U_\beta=\setof{\bbu_\xi}{\xi<\beta}$. We show that 
$X_\beta=[U_\beta]_{\ell_2(\kappa)}$ is non-pathological. 

Let $\mapping{f}{E\cap\beta}{\beta}$ be as in \xitemof{ads-1}. 
For each $\alpha\in E\cap\beta$, let
$B_\alpha=(A_\alpha\setminus f(\alpha))\cup\ssetof{\alpha}$. Then $B_\alpha$,
$\alpha\in E\cap\beta$ are pairwise disjoint. Let
$C=\beta\setminus(\bigcup_{\alpha\in E\cap\beta}B_\alpha)$. 

Note that 
\begin{xitemize}
\xitem[ads-2] $\bbu'_\alpha=\bbu_\alpha
  -\sum_{\xi\in A_\alpha\cap f(\alpha)}\bbu_\alpha(\xi)\bbe^\kappa_\xi$
\end{xitemize}
is an element of $X$ and $\supp(\bbu'_\alpha)=B_\alpha$. It follows that $X_\beta$ is 
the orthogonal sum of the sub-inner-product-spaces $X\downarrow C$, $X\downarrow B_\alpha$,
$\alpha\in E\cap\beta$. In particular, we have
$X_\beta=\oplusbar_{\alpha\in E\cap\beta}X\downarrow B_\alpha\ \oplus\  X\downarrow C$. 
Note that from this it follows that $X_\beta=X\downarrow\beta$. 

Now $X\downarrow B_\alpha$, $\alpha\in E\cap\beta$ 
are non-pathological since they are separable. Let $U_\alpha$ be an orthonormal 
basis of $X\downarrow B_\alpha$ for $\alpha\in E\cap\beta$. Also $X\downarrow  C$ is 
non-pathological with the orthonormal basis
$\setof{\bbe^\kappa_\alpha}{\alpha\in C}$. Thus
$\bigcup_{\alpha\in E\cap\beta}U_\alpha\cup\setof{\bbe^\kappa_\alpha}{\alpha\in C}$
is an orthonormal basis of $X_\beta$. 

The same argument shows that $X\downarrow S$ is non-pathological for any 
bounded subset $S$ of $\kappa$ closed \wrt\ 
the sequence $\seqof{A_\alpha}{\alpha\in E}$ (that is,  
$A_\alpha\subseteq S$ for all $\alpha\in E\cap S$).  Note that, for all regular 
$\lambda<\kappa$ there are club many such $S$ of cardinality $\lambda$. 
\qedofThm\qedskip

Under the consistency strength of certain very large cardinals we obtain reflection 
theorems for pathology of pre-Hilbert spaces.

\begin{Thm}
  \label{T-refl-2}
  Suppose that $\kappa$ is a supercompact cardinal. Then for any pathological 
  pre-Hilbert space $X$, there are stationarily many pathological 
  sub-inner-product-spaces $Y$ of $X$ of size $<\kappa$. 
\end{Thm}
\prf Suppose that $X$ is a pathological pre-Hilbert space of size $\lambda$. 
We may assume that the underlying set of $X$ is $\lambda$. 
If $\lambda<\kappa$ then the statement of the theorem is trivial. So we assume that
$\lambda\geq\kappa$. Let $\calC\subseteq[\lambda]^{<\kappa}$ be a club set.
Let $\elembed{j}{V}{M}$ be an elementary embedding with $\crit(j)=\kappa$, 
$j(\kappa)>\lambda$ and $\fnsp{\lambda}{M}\subseteq M$. Then we have
$j\imageof X\in M$ and $j\imageof X\in j(\calC)$: the latter is because 
$M\models\mbox{``}j(\calC)$ is a club subset of $[j(\lambda)]^{<j(\kappa)}$'' and 
$\calD=\setof{j(Y)}{Y\in\calC}\subseteq j(\calC)$ is of 
cardinality $<j(\kappa)$ with $\bigcup\calD=j\imageof X$.  
We have $V\models j\imageof X\cong X$ and hence
$V\models\mbox{``}j\imageof X\mbox{ is pathological''}$. It 
follows that $M\models\mbox{``}j\imageof X\mbox{ is pathological''}$.
Putting these facts together, we obtain
\begin{xitemize}
\xitem[refl-0] 
  $M\models$``\
  \parbox[t]{\textwidth}{$j\imageof X$ is a sub-inner-product space of $j(X)$,
  $j\imageof X\in j(\calC)$ and\\
  $j\imageof X$ is pathological''.} 
\end{xitemize}

Thus, 
\begin{xitemize}
\xitem[refl-1] 
  $M\models$``there is a pathological sub-inner-product-space $Y$ of $j(X)$ with
  $Y\in j(\calC)$. 
\end{xitemize}
By elementarity if follows 
\begin{xitemize}
\xitem[refl-2] 
  $V\models$``there is a pathological sub-inner-product-space $Y$ of $X$ with $Y\in \calC$.
\qedofThm\qedskip
\end{xitemize}

\begin{Thm}
  \label{T-refl-3}
  Suppose that $X$ is a pathological pre-Hilbert space and $X$ is a dense 
  sub-inner-product-space of $\ell_2(S)$ for some infinite set $S$. Then for any 
  ccc \po\ $\poP$ we have
  $\forces{\poP}{[X]_{\ell_2(S)}\mbox{ is pathological\/}}$. 
\end{Thm}
\prf Suppose that $X$ is a pre-Hilbert space and there is a ccc \po\ $\poP$ \st\ 
\begin{xitemize}
\xitem[refl-3] 
  $\forces{\poP}{[X]_{\ell_2(S)}\mbox{ is non-pathological\/}}$. 
\end{xitemize}
We show that $X$ is then non-pathological.

By \Thmof{L-char-2} and the Maximal Principle there is 
a $\poP$-name $\utilde{\calP}$ of partition of $S$ into countable sets \st\ 
\begin{xitemize}
\xitem[refl-4] 
  $\forces{\poP}{[X]_{\ell_2(S)}
  =\oplusbar_{P\in\utilde{\calP}}[X]_{\ell_2(S)}\downarrow P}$.
\end{xitemize}
\begin{Claim}
  \label{T-refl-3-C-1}
  There is a partition $\calP'$ of $\kappa$ into countable sets \st, for each
  $P\in\calP'$, we have
  $\forces{\poP}{P\mbox{ is a countable union of elements of }\utilde{\calP}}$. 
\end{Claim}
\prfofClaim
Let $\sim$ be the transitive closure of the relation
\begin{xitemize}
\xitem[refl-5] 
  $\sim_0=\setof{\pairof{s,t}\in S}{
  \begin{array}[t]{@{}l}
    \mbox{there is }p\in\poP\mbox{ \st}\\[\jot]
    p\forces{\poP}{s\mbox{ and }t\mbox{ belong to the same 
        set }\in\utilde{\poP}}.}
  \end{array}
$ 
\end{xitemize}

By the ccc of $\poP$, $Q_s=\setof{t\in S}{s\sim_0 t}$ is countable for all
$s\in S$. Hence all equivalence classes of $\sim$ are also countable. 

Let $\calP'$ be the partition of $S$ into equivalence classes of $\sim$. 

Let $P\in\calP'$. We show that
$\forces{\calP}{P\xmbox{ is a union of elements of }\utilde{\calP}}$. Let $G$ be 
an arbitrary $(V,\poP)$-generic set In $V[G]$ suppose that $s\in P$, $s\in Q$ for 
some $Q\in\utilde{\calP}^G$ and $t\in Q$. Then there is some $p\in G$ \st\
$p\forces{\poP}{s\mbox{ and }t\mbox{ are in the same element of }\utilde{\calP}}$.  
It follows that $s\sim_0 t$ and $t\in P$. 
\qedofClaim
\qedskip

\begin{Claim}
  \label{T-refl-3-C-2}
  For $P\in\calP'$ we have $X=(X\downarrow P)\oplus(X\downarrow (S\setminus P))$.
\end{Claim}
\prfofClaim Let $G$ be a $(V,\poP)$-generic set. In $V[G]$, we have
$[X]_{\ell_2(S)}
=([X]_{\ell_2(S)}\downarrow P) \oplus ([X]_{\ell_2(S)}\downarrow (S\setminus P))$ by 
\Claimof{T-refl-3-C-1}. Hence, in $V$, we have 
$X=(X\downarrow P) \oplus (X\downarrow (S\setminus P))$
\qedofClaim
\qedskip

It follows from \Claimabove\ that 
$X=\oplusbar_{P\in\calP'}X\downarrow P$. Thus, by \Thmof{T-char}, $X$ is has an 
orthonormal basis. \qedofThm
\qedskip

The Cohen forcing $\Fn(\kappa,2)$ in the following theorem can be replaced by may 
other c.c.c.\ forcing notions which can be seen as iterations with certain 
coherence (see Dow, Tall Weiss \cite{dow-etal}).

\begin{Thm}
  Assume that $\kappa$ is a supercompact cardinal and let $\poP=\Fn(\kappa,2)$. 
  Then we have
  \begin{xitemize}
  \xitem[refl-6] 
    $\forces{\poP}{
    \begin{array}[t]{@{}l}
      \mbox{for every pathological pre-Hilbert space }X\mbox{ which is a dense sub-inner-}\\
      \mbox{product space of }\ell_2(\lambda)\mbox{ for some infinite }\lambda,
      \mbox{ there are stationarily many}\\
      S\in[\lambda]^{<2^{\aleph_0}}\mbox{ \st\ }
      X\downarrow S\mbox{ is pathological\/}}.
    \end{array}$
  \end{xitemize}
\end{Thm}
\prf Let $G$ be a $(V,\poP)$-generic filter. Working in $V[G]$, let $X$ be a 
pre-Hilbert space which is a dense sub-inner-product-space of $\ell_2(\lambda)$. 
If $\lambda<\kappa=\left(2^{\aleph_0}\right)^{V[G]}$ then the assertion is 
trivial. Thus we assume $\lambda\geq\kappa$. Let $\calC\subseteq[\lambda]^{<\kappa}$ be a 
club set. It is enough to show that there is some $S\in\calC$ \st\ $X\downarrow S$ 
is pathological. 

Back in $V$, let $\elembed{j}{V}{M}$ be a $\lambda$-supercompact embedding. That 
is, the elementary embedding $j$ is \st\ $M\subseteq V$ is a transitive class
$\crit(j)=\kappa$, $j(\kappa)>\lambda$ and $\fnsp{\lambda}{M}\subseteq M$. 

Let $\poP^*=\Fn(j(\kappa),2)=j(\poP)$ and let $G^*$ be a $(V,\poP^*)$-generic 
filter with $G^*\supseteq G$. 

Let $\elembed{j^*}{V[G]}{M[G^*]}$ be the extension of $j$ defined by
\begin{xitemize}
\xitem[refl-7] 
  $j^*([\utilde{a}]^G)=[j(\utilde{a})]^{G^*}$ 
\end{xitemize}
for each $\poP$-name $\utilde{a}$. It is easy to check that $j^*$ is 
well-defined and 
\begin{xitemize}
\xitem[refl-8] 
  $\left(\fnsp{\lambda}{M[G^*]}\right)^{V[G^*]}\subseteq M[G^*]$. 
\end{xitemize}
It follows that $j^*\imageof X\in M[G^*]$ and
$\supp(j^*\imageof X)=j\imageof\lambda\in j^*(\calC)$.  
Since $V[G^*]$ is a c.c.c.\ extension of $V[G]$, 
by 
\Lemmaof{T-refl-3}, we have 
\begin{xitemize}
\xitem[] 
  $V[G^*]\models\mbox{``} [j\imageof X]_{\ell_2(j(\lambda))}\mbox{ is pathological''}$. 
\end{xitemize}

It follows that 
\begin{xitemize}
\xitem[refl-9] 
  $M[G^*]\models\mbox{``} [j\imageof X]_{\ell_2(j(\lambda))}\mbox{ is pathological''}$ 
\end{xitemize}
by the same argument as right after \xitemof{refl-0}. 

Thus we have
\begin{xitemize}
\xitem[refl-10] 
  $M[G^*]\models$``there is $S\in j^*(\calC)$ \st\ $j(X)\downarrow S$ is pathological''.
\end{xitemize}
By elementarity it follows that 
\begin{xitemize}
\xitem[refl-11] 
  $V[G]\models$``there is $S\in\calC$ \st\ $X\downarrow S$ is pathological''.
\end{xitemize}
\qedofThm
\section{A Singular Compactness Theorem}
\label{singular}
The proof of the following theorem follows closely the proof of 
Shelah's Singular Compactness Theorem given in Hodges \cite{hodges}. A similar 
Singular Compactness Theorem in the context of (non-)freeness of modules is given 
in Eklof \cite{eklof}. 

\begin{Thm}
  \label{Th-SC}
  Suppose that $\lambda$ is a singular cardinal and $X$ is a pre-Hilbert space 
  which is a dense sub-inner-product-space of $\ell_2(\lambda)$. If $X$ is 
  pathological then there is a cardinal $\lambda'<\lambda$ \st\ 
  \begin{xitemize}
  \xitem[sc-1] 
    $\setof{u\in[\lambda]^{\kappa^+}}{X\downarrow u\mbox{ is a pathological 
      pre-Hilbert space}}$
  \end{xitemize}
  is stationary in $[\lambda]^{\kappa^+}$ for all $\lambda'\leq\kappa<\lambda$. 
\end{Thm}
In the following we shall prove the contraposition of 
the statement of the theorem:

\begin{list}{}{\setlength{\leftmargin}{0pt}}
  \item {\bf \bfThmof{Th-SC}${}^{\mbox{\bf*}}$}
  \it 
  For any singular $\lambda$ and any pre-Hilbert space $X$ which is a dense 
  sub-inner-product-space of $\ell_2(\lambda)$, if 
  \item
  \begin{xitemize}
  \xitem[sc-2] 
    $\calN^X_{\kappa}=\setof{u\in[\lambda]^{\kappa^+}}{X\downarrow u
  \mbox{ is a non-pathological pre-Hilbert space}}$  contains a club in
    $[\lambda]^{\kappa^+}$ for cofinally many $\kappa<\lambda$, 
  \end{xitemize}
  then $X$ is non-pathological.
\end{list}


For a dense sub-inner-product space $X$ of $\ell_2(\lambda)$ and $v$,
$v'\subseteq \lambda$, we write $u'\osmmdin{X}u$ if 
$u\subseteq u'$, $X\downarrow u$ and $X\downarrow u'$ are dense in
$\ell_2(\lambda)\downarrow u$ and $\ell_2(\lambda)\downarrow u'$ respectively; 
and 
$X\downarrow u$ is an orthogonal direct summand of $X\downarrow u'$, i.e.\ if 
$X\downarrow u'=(X\downarrow u)\oplus (X\downarrow (u'\setminus u))$, see 
\Lemmaof{summand}.  

For a cardinal $\kappa$, the $\kappa$-Shelah game over 
$X\subseteq\ell_2(\lambda)$ (notation $\calG_\kappa(X)$) is the game whose 
matches $\calM$ are $\omega$-sequences of  
moves by Players I and II
\begin{xitemize}
\item[] $\calM:\ \ \ 
  \begin{array}{lllll}
    {\rm I}\quad &u_0 &u_1 &u_2 &\cdots\\
    {\rm II}\quad &v_0 &v_1 &v_2 &\cdots\\
  \end{array}
  $
\end{xitemize}
where $u_i$, $v_i\in[\lambda]^{\kappa}$ for $i\in\omega$ and
$u_0\subseteq v_0\subseteq u_1\subseteq v_1\subseteq u_2\subseteq v_2\subseteq\cdots$.

Player II wins if $X\downarrow v_i$ is non-pathological and 
$v_{i+1}\osmmdin{X}v_i$ for all $i\in\omega$. 

Note that, if Player II wins in a match $\calM$ with the moves
$u_0\subseteq v_0\subseteq u_1\subseteq v_1\subseteq u_2\subseteq v_2\subseteq\cdots$,
then $X\downarrow w$ for $w=\bigcup_{i\in\omega}u_i=\bigcup_{i\in\omega}v_i$ is 
non-pathological.

\begin{Lemma}
  \label{L-sc-0}
  $\kappa$-Shelah game over $X\subseteq\ell_2(\lambda)$ is determined for 
  regular $\kappa$.  
\end{Lemma}
\prf Since the game is open for Player I, the proof of Gale-Stewart Theorem 
applies (see e.g.\ Kanamori \cite{kanamori} or Hodges \cite{hodges}). \qedofLemma
\qedskip

\begin{Lemma}
  \label{L-sc-1}
  Suppose that $X$ is 
  a dense sub-inner-product-space of $\ell_2(\lambda)$ for a cardinal
  $\lambda$.  For a cardinal $\kappa<\lambda$,  if $\calN^X_{\kappa}$ contains a 
  club subset of $[\lambda]^{\kappa^+}$,  
  then Player II has a winning strategy in $\calG_\kappa(X)$. 
\end{Lemma}
\prf By \Lemmaof{L-sc-0}, it is enough to show that the Player I does not have a 
winning strategy. 

Suppose that $\sigma$ is a strategy for Player I. We show that it is not 
winning. 

Let $\calC\subseteq\calM^X_\kappa$ be club in $[\lambda]^{\kappa^+}$. 

Let $\chi$ be a sufficiently large regular cardinal and let 
$\seqof{M_\alpha}{\alpha<\kappa^+}$ be a continuously increasing chain of 
elementary submodels of $\calH(\theta)$ \st\ 
\begin{xitemize}
\xitem[sc-2-0] $\sigma$, $X$, $\lambda$, $\kappa$, $\calC\ctenten$\ $\in M_0$;
\xitem[sc-3] $\cardof{M_\alpha}=\kappa$ and $M_\alpha\in M_\alpha+1$ for all $\alpha<\kappa^+$;
\xitem[sc-3-0] $\alpha\subseteq M_\alpha$ for all $\alpha<\kappa^+$
\xitem[sc-4] For any finite subsequence $\calG_0$ of
  $\seqof{M_\beta}{\beta\leq\alpha}$, if $\calG_0$ is the moves of Player II in an 
  initial segment $\calM_0$ of a match in $\calG_\kappa(X)$ where the Player I 
  has played according to $\sigma$ and the last member of $\calG_0$ is the last 
  move in $\calM_0$, then $\sigma(\calM_0)\in M_{\alpha+1}$ and
  $\sigma(\calM_0)\subseteq M_{\alpha+1}$. 
\end{xitemize}

Let $M=\bigcup_{\alpha<\kappa^+}M_\alpha$. By \xitemof{sc-2-0}, \xitemof{sc-3} 
and \xitemof{sc-3-0}, we have 
\begin{xitemize}
\xitem[sc-5] $\lambda\cap M\in\calC$.
\end{xitemize}
By \Thmof{T-char}, there is a partition $\calP$ of $\lambda\cap M$ into countable 
sets \st\ $X\downarrow(\lambda\cap M)=\oplusbar_{A\in\calP}X\downarrow A$. 
Let
$C=\setof{\alpha<\kappa^+}{\lambda\cap M_\alpha\xmbox{ is a union of some elements of }
\calP}$. Then $C$ is a club set $\subseteq\kappa^+$,
\begin{xitemize}
\xitem[sc-5-0] 
  $M\downarrow(\lambda\cap M_\alpha)$ is non-pathological for all $\alpha\in C$ and 
\xitem[sc-6] 
  $(\lambda\cap M_\alpha)\osmmdin{X}(\lambda\cap M)$ for every $\alpha\in C$. 
\end{xitemize}

Let $\alpha_i$, $i\in\omega$ be the first $\omega$ elements of $C$ and
$v_i=\lambda\cap M_{\alpha_i}$ for $i\in\omega$. 
By \xitemof{sc-4}, there is a match $\calM$ in $\calG_\kappa(X)$ in which Player 
I has chosen his moves according to $\sigma$ and $\seqof{v_i}{i\in\omega}$ is 
the moves of Player II. Player II wins in this match $\calM$ by \xitemof{sc-6}. 
This shows that $\sigma$ is not a winning strategy of Player I.\qedofLemma
\qedskip

\noindent
{\bf Proof of \bfThmof{Th-SC}${}^{\mbox{\bf*}}$:}\ \ Suppose that $X$ and $\lambda$ are as in 
\Thmof{Th-SC}${}^{\mbox{*}}$.  Let $\delta=\cf(\lambda)$ and let 
$\seqof{\lambda_\xi}{\xi<\delta}$ be a continuously increasing sequence of 
cardinals below $\lambda$ \st\
\begin{xitemize}
\xitem[sc-7] $\delta<\lambda_0$;
\xitem[sc-8] $\calN^X_{\lambda_\xi}$ (defined in \xitemof{sc-2}) contains a club subset
  $\subseteq[\lambda]^{(\lambda_\xi)^+}$ for all successor $\xi<\delta$. 
\end{xitemize}
The condition \xitemof{sc-8} is possible by our assumption \xitemof{sc-7}.

In the following, we construct $u^i_\xi$, $\tilde{u}^i_\xi$, $v^i_\xi$ for 
$\xi<\delta$ and $i\in\omega$ \st\ 

\begin{xitemize}
\xitem[sc-9]
  $\lambda_\xi=u^0_\xi\subseteq\tilde{u}^0_\xi\subseteq v^0_\xi\subseteq
  u^1_\xi\subseteq\tilde{u}^1_\xi\subseteq v^1_\xi\subseteq
  u^2_\xi\subseteq\tilde{u}^2_\xi\subseteq v^2_\xi\subseteq\cdots$
\end{xitemize}
and, 
letting $w_\xi=\bigcup_{i\in\omega}u^i_\xi=\bigcup_{i\in\omega}\tilde{u}^i_\xi
=\bigcup_{i\in\omega}v^i_\xi$, we have
\begin{xitemize}
\xitem[sc-10] $\seqof{w_\xi}{\xi\in\delta}$ is a filtration of $\lambda$;
\xitem[sc-11] $X\downarrow w_\xi$ is non-pathological for all $\xi\in\delta$;
\xitem[sc-12] $w_\eta\osmmdin{X}w_\xi$ for all $\xi<\eta<\delta$.
\end{xitemize}
From \xitemof{sc-10}, \xitemof{sc-11} and \xitemof{sc-12}, it follows immediately 
that $X$ is non-pathological.

For the construction of $u^i_\xi$, $\tilde{u}^i_\xi$, $v^i_\xi$ for 
$\xi<\delta$ and $i\in\omega$, we fix winning strategies $\sigma_\xi$ for Player 
II in $\calG_{\lambda_\xi}(X)$ for all successor $\xi<\delta$. We have such strategies by 
\xitemof{sc-8} and \Lemmaof{L-sc-1}. 

The following describes the inductive construction:
\begin{xitemize}
\xitem[sc-13]
  $\cardof{u^i_\xi}=\cardof{\tilde{u}^i_\xi}=\cardof{v^i_\xi}=\lambda_\xi$ for 
  all $\xi<\delta$;
\xitem[sc-14] The sequence $\tilde{u}^0_\xi$, $v^0_\xi$, $\tilde{u}^1_\xi$,
  $v^1_\xi$, $\tilde{u}^2_\xi$, $v^2_\xi$\ctenten\ is a match in
  $\calG_{\lambda_\xi}(X)$ in which Player II has played according 
  to $\sigma_\xi$ for all successor $\xi<\delta$ (\xitemof{sc-11} for all 
  successor $\xi<\delta$ follows from this); 
\xitem[sc-15]  When 
$\seqof{u^k_\xi}{k\leq i,\,\xi<\delta}$,
  $\seqof{\tilde{u}^j_\xi}{j<i,\,\xi<\delta}$ and
  $\seqof{v^j_\xi}{j<i,\,\xi<\delta}$ have been chosen (according to all the 
  conditions described here) for 
  an $i\in\omega$ then $\tilde{u}^i_\xi$ for each $\xi<\delta$ is \st\
  $\tilde{u}^i_\xi\supseteq\bigcup_{\eta\leq\xi}u^i_\eta$ holds
  (note that $\cardof{\bigcup_{\eta\leq\xi}u^i_\eta}=\lambda_\xi$ by 
  \xitemof{sc-13}. This condition guarantees that the sequence
  $\seqof{w_\xi}{\xi<\delta}$ is going to be increasing); 
\end{xitemize}

For each successor $\xi<\delta$ and $i\in\omega$, if $v^i_\xi$ has been chosen according to 
the conditions described here, $X\downarrow v^i_\xi$ is non-pathological by 
\xitemof{sc-14}. Thus we can find a partition $\calP^i_\xi$ of $v^i_\xi$ into 
countable sets \st\ $X\downarrow v^i_\xi=\oplusbar_{A\in\calP^i_\xi}X\downarrow A$
by \Thmof{T-char}. If $i>0$ then we may choose $\calP^i_\xi$ \st\
$\calP^{i-1}_\xi\subseteq \calP^i_\xi$ (this is possible since
$v^i_\xi\osmmdin{X}v^{i-1}_\xi$ by \xitemof{sc-14}). 

\begin{xitemize}
\xitem[sc-16] (a continuation of \xitemof{sc-15}) When 
$\seqof{u^k_\xi}{k\leq i,\,\xi<\delta}$,
  $\seqof{\tilde{u}^j_\xi}{j<i,\,\xi<\delta}$ and
  $\seqof{v^j_\xi}{j<i,\,\xi<\delta}$ have been chosen (according to all the 
  conditions described here) for 
  an $i\in\omega$ then we choose $\tilde{u}^i_\xi$ also \st\ 
$\tilde{u}^i_\xi\cap v^k_{\xi+1}$ is a union of some elements of $\calP^k_\xi$ for 
  all $k<i$ for all (not necessarily successor) $\xi<\delta$
  (this makes $w_{\xi+1}\osmmdin{X}w_\xi$ for all $\xi<\delta$);
\end{xitemize}

For each $\xi<\delta$ and $i\in\omega$, when $v^i_\xi$ has been chosen, we enumerate it as
$v^i_\xi=\setof{\beta_{i,\xi,\eta}}{\eta<\lambda_\xi}$. 
\begin{xitemize}
\xitem[sc-17] When $\seqof{u^j_\xi}{j< i,\,\xi<\delta}$,
  $\seqof{\tilde{u}^j_\xi}{j<i,\,\xi<\delta}$ and
  $\seqof{v^j_\xi}{j<i,\,\xi<\delta}$ have been chosen (according to all the 
  conditions described here) for 
  an $i\in\omega$ then we let
  $u^{i+1}_\xi=\setof{\beta_{i,\xi,\eta}}{\xi<\delta,\,\eta<\lambda_\xi}\cup v^i_\xi$ \\
  (this is possible since the set on the right side of the inequality has size
  $\leq\lambda_\xi$. This condition makes the 
  sequence $\seqof{w_\xi}{\xi<\delta}$ continuous). 
\end{xitemize}

To see that \xitemof{sc-17} makes the sequence $\seqof{w_\xi}{\xi<\delta}$ 
continuous, suppose that $\nu\in w_\gamma$ for a limit $\gamma<\delta$. 
Then there is $i^*\in\omega$ \st\ $\nu\in v^{i^*}_\gamma$. Hence there is 
$\eta^*<\lambda_\gamma$ \st\ $\nu=\beta_{i^*,\gamma,\eta^*}$. Let $\xi<\gamma$ 
be \st\ $\eta^*<\lambda_\xi$. Then by \xitemof{sc-17} we have
$\nu=\beta_{i^*,\gamma,\eta^*}\in u^{i^*+1}_\xi\subseteq w_\xi$. 

As noted above, the choice of  
$u^i_\xi$, $\tilde{u}^i_\xi$, $v^i_\xi$ for 
$\xi<\delta$ and $i\in\omega$ with 
\xitemof{sc-9}, \xitemof{sc-13} $\sim$ \xitemof{sc-17} makes
$\seqof{w_\xi}{\xi<\delta}$ satisfy the conditions 
\xitemof{sc-10}, \xitemof{sc-11} for all successor $\xi<\delta$ and 
\xitemof{sc-12} for all $\xi<\delta$ and $\eta=\xi+1$. 

By the continuity of 
$\seqof{w_\xi}{\xi<\delta}$ we can then prove inductively that \xitemof{sc-11} and 
\xitemof{sc-12} hold for all $\xi<\eta<\delta$. \qedof{\Thmof{Th-SC}${}^{\mbox{*}}$}

\section{Reflection of pathology and Fodor-type Reflection 
  Principle}
\label{FRP}
In this section we prove the following theorem which gives characterizations of 
\FRP\ in terms of pathology of pre-Hilbert spaces.

\begin{Thm}
  \label{main-thm-1} Each of the following assertions is equivalent to \FRP:
  \begin{xitemize}
  \xitem[chr-23] For any regular $\kappa>\omega_1$ and any dense sub-inner-product-space 
    $X$ of $\ell_2(\kappa)$, if $X$ is pathological then 
    \begin{xitemize}
    \item[] $S_X=\setof{\alpha<\kappa}{X\downarrow \alpha\mbox{ is pathological\/}}$
    \end{xitemize}
    is stationary in $\kappa$. 
  \xitem[chr-24] For any regular $\kappa>\omega_1$ and any dense sub-inner-product-space 
    $X$ of $\ell_2(\kappa)$, if $X$ is pathological then 
    \begin{xitemize}
    \item[] $S^{\aleph_1}_X=\setof{U\in[\kappa]^{\aleph_1}}{
      X\downarrow U\mbox{ is pathological\/}}$
    \end{xitemize}
    is stationary in $[\kappa]^{\aleph_1}$. 
  \end{xitemize}
\end{Thm}

First let us review some facts around the reflection principle \FRP\ needed 
for the proof of \Thmof{main-thm-1}. 

One of the combinatorial statements equivalent to \FRP\ we are going to use below is as 
follows:
\begin{xitemize}
\item[(\FRP)] For any regular $\kappa>\omega_1$, any stationary 
  $E\subseteq E^\omega_\kappa$ and any mapping $\mapping{g}{E}{[\kappa]^{\aleph_0}}$, 
  there is $\alpha^*\in E^{\omega_1}_\kappa$ \st\ 
  \begin{xitemize}
    \xitem[frp-0] 
    $\alpha^*$ is closed \wrt\ $g$ 
    (that is, $g(\alpha)\subseteq\alpha^*$ for all $\alpha\in E\cap\alpha^*$) and, for 
    any $I\in[\alpha^*]^{\aleph_1}$ closed \wrt\ $g$, closed in $\alpha^*$ \wrt\ 
    the order topology and with $\sup(I)=\alpha^*$, if 
    $\seqof{I_\alpha}{\alpha<\omega_1}$ is a filtration of $I$ then
    $\sup(I_\alpha)\in E$ and 
    $g(\sup(I_\alpha))\cap\sup(I_\alpha)\subseteq I_\alpha$ hold for stationarily many
    $\alpha<\omega_1$  
  \end{xitemize}
\end{xitemize}
(see Fuchino, Sakai, Soukup \cite{more}). 

\FRP\ was invented by Lajos Soukup and the author in 2008 and then published in 
Fuchino Juha\'asz, Soukup,  
  Szentmikl\'ossy, Usuba \cite{fuchino-juhasz-etal} by a formulation slightly 
different from the one given above. 
In Fuchino, Sakai, Soukup \cite{more} it is proved that \FRP\ is equivalent to 
the statement that  
$\ADS^-(\kappa)$ fails for all regular $\kappa>\omega_1$. This characterization 
of \FRP\ is used to show the equivalence of \FRP\ to many mathematical reflection 
statements
in Fuchino \cite{balogh}, Fuchino, Sakai, Soukup \cite{more}, Fuchino, Rinot \cite{rinot}. 
One of the typical mathematical assertion equivalent with \FRP\ is:
\begin{xitemize}
\item[] For every non-metrizable countably compact topological space $X$ there is 
  a non-metrizable subspace of $X$ of cardinality $\leq\aleph_1$ (see \cite{more}). 
\end{xitemize}
Our present result adds another couple of mathematical 
reflection statements to the long list of the statements equivalent to \FRP. 

For the proof of \Thmof{main-thm-1} we need the following easy observations:

\begin{Lemma}{\rm (cf. Lemma 6.1 in \cite{fuchino-juhasz-etal})}
\label{frp-2}
  Suppose that $\kappa$ is a regular cardinal $>\aleph_1$, $C\subseteq\kappa$ 
  club, $E\subseteq C$ 
  stationary and $a_\eta\in[\kappa]^{\aleph_0}$ for $\eta\in E$. Then there is a 
  stationary $E'\subseteq E^\omega_\kappa\cap C$ and a mapping 
  $\mapping{\overline{\eta}}{E'}{E}$; $\xi\mapsto\eta_\xi$ 
  \st, for all $\xi\in E'$, we 
  have $\xi\leq\eta_\xi$ and $a_{\eta_\xi}\cap\xi=a_{\eta_\xi}\cap\eta_\xi$. 
\end{Lemma}
\prf 
We prove the Lemma in the following two cases:

\noindent
{\bf Case I.} $E\cap E^\omega_\kappa$ is stationary. 

Then 
$E'=E\cap E^\omega_\kappa$ with $\bar{\eta}=\id_{E'}$ is as 
desired. \smallskip

\noindent
{\bf Case II.} $E\cap E^\omega_\kappa$ is non-stationary. 
Then $E''=E\setminus E^\omega_\kappa$ is stationary. For each $\eta\in E''$ we have
$\sup(a_\eta\cap\eta)<\eta$. By Fodor's Lemma there are $\eta_0<\kappa$ and stationary
$E'''\subseteq E''$ \st\ $\sup(a_\eta\cap\eta)\leq\eta_0$ for all $\eta\in E'''$. 

Let $E'=(E^\omega_\kappa\cap C)\setminus \eta_0$ and, for each $\xi\in E'$, let
$\eta_\xi=\min(E'''\setminus \xi)$. Then this $E'$ with
$\mapping{\bar{\eta}}{E'}{E}$; $\xi\mapsto\eta_\xi$ is as desired. 
\qedofLemma

\begin{Lemma}
  \label{frp-2-0}
  Suppose that $\kappa$ and $\lambda$ are regular cardinals with 
  $\aleph_0<\kappa<\lambda$ and $A$ a set of size $\geq\lambda$. 
  If $S\subseteq[A]^{<\lambda}$ is stationary in $[A]^{<\lambda}$ and 
  $U_s\subseteq[s]^{<\kappa}$ is stationary for all $s\in S$, then 
  $\bigcup_{s\in S}U_s$ is stationary in $[A]^{<\kappa}$. 
\end{Lemma}
\prf
Suppose that $C\subseteq[A]^{<\kappa}$ is club. Then there is 
$\mapping{f}{[A]^{<\aleph_0}}{[A]^{<\kappa}}$ \st\ 
\begin{xitemize}
\xitem[frp-2-1]
  $C_f=\setof{x\in[A]^{<\kappa}}{x\mbox{ is closed \wrt\ }f}\subseteq C$ 
\end{xitemize}
(see e.g. Lemma 8.26 in Jech \cite{millenium-book}). 
Note that then
\begin{xitemize}
\xitem[frp-2-2] 
  $C^{<\lambda}_f=\setof{y\in[A]^{<\lambda}}{y\mbox{ is closed \wrt\ }f}$
\end{xitemize}
is a club $\subseteq[A]^{<\lambda}$. 

Since $S$ is a stationary subset of $[A]^{<\lambda}$, there is
$s^*\in S\cap C^{<\lambda}_f$. Now  
$C_f\cap[s^*]^{<\kappa}$ is a club in $[s^*]^{<\kappa}$ and $U_{s^*}$ is 
stationary in $[s^*]^{<\kappa}$. 

Thus there is 
$u^*\in U_{s^*}\cap (C_f\cap[s^*]^{<\kappa})\subseteq (\bigcup_{s\in S}U_s)\cap C_f
\subseteq (\bigcup_{s\in S}U_s)\cap C$. \\
\qedofLemma
\qedskip

\noindent
{\bf Proof of \bfThmof{main-thm-1}:}  First we show that \FRP\ implies \xitemof{chr-23}.

Assume that \FRP\ holds. 
Suppose that $X$ is a dense 
sub-inner-product-space of $\ell_2(\kappa)$ for a regular cardinal 
$\kappa>\aleph_1$. We assume that $S_X$ (in \xitemof{chr-23}) is non-stationary 
and drive a contradiction. 

By the assumption there is a club set $C\subseteq\kappa$ \st\ $X\downarrow\alpha$ 
is non-pathological for all $\alpha\in C$. By \Lemmaof{L-dense-1} we may assume 
that $X\downarrow\alpha$ is dense in $\ell_2(\kappa)\downarrow\alpha$ for all
$\alpha\in C$. 

Since $X$ is pathological, 
\begin{xitemize}
\xitem[frp-3] 
  $E=\setof{\alpha\in C}{X\downarrow\alpha
  \mbox{ is not an orthogonal direct summand of }X}$
\end{xitemize}
is stationary. By \Lemmaof{summand}, there is $\bba_\alpha\in X$ \st\
$\bba_\alpha\downarrow\alpha\not\in X\downarrow\alpha$ for all $\alpha\in E$. 
Let $A_\alpha=\supp(\bba_\alpha)$ for $\alpha\in E$. By 
\Lemmaof{frp-2}, we may assume \wolog\ that $E\subseteq C\cap E^\omega_\kappa$. 

By \FRP, there is $\alpha^*\in E^{\omega_1}_\kappa$ \st\ \xitemof{frp-0} holds 
for $\mapping{g}{E}{[\kappa]^{\aleph_0}}$; $\alpha\mapsto A_\alpha$. 

Now since $E\cap\alpha^*$ is unbounded in $\alpha^*$, we have 
$\alpha^*\in C$. Thus $X\downarrow\alpha^*$ is non-pathological. Hence by 
\Thmof{T-char} there are club many $I\in[\alpha^*]^{\aleph_1}$ \st\ 
$X\downarrow I$ is non-pathological and $X\downarrow I$ is dense in
$\ell_2(\kappa)\downarrow I$. It follows that there is $I^*\in[\alpha^*]^{\aleph_1}$ 
\st\ 
\begin{xitemize}
\xitem[frp-4] $I^*$ is closed \wrt\ $g$ and closed in $\alpha^*$ \wrt\ the order topology;
\xitem[frp-4-0] $X\downarrow I^*$ is non-pathological;
\xitem[frp-4-1] $X\downarrow I^*$ is dense in $\ell_2(\kappa)\downarrow I^*$ and 
\xitem[frp-5] $\sup(I^*)=\alpha^*$.
\end{xitemize}

Let $\seqof{I_\alpha}{\alpha<\omega_1}$ be a filtration of $I^*$ \st\
$X\downarrow I_\alpha$ is dense in $\ell_2(\kappa)\downarrow I_\alpha$. 

By \xitemof{frp-0}, 
\begin{xitemize}
\xitem[frp-6] 
  $E_0=\setof{\alpha\in\omega_1}{\sup(I_\alpha)\in E,\,
  A_{\sup(I_\alpha)}\cap\sup(I_\alpha)\subseteq I_\alpha}$
\end{xitemize}
is stationary. By \Lemmaof{filtration} and \Lemmaof{summand} this is a 
contradiction to \xitemof{frp-4-0}. 
This proves that \FRP\ implies \xitemof{chr-23}. 

Since \FRP\ is equivalent to the 
global negation of $\ADS^-(\kappa)$. \Thmof{T-refl-1} implies the converse. 

For the equivalence of \FRP\ and \xitemof{chr-24}, it is enough by virtue of the 
second part of \Thmof{T-refl-1} to show that 
\xitemof{chr-23} implies \xitemof{chr-24}.

Assume that \xitemof{chr-23} holds. We prove that \xitemof{chr-24} holds for all 
uncountable $\kappa$ by induction on $\kappa$: if $\kappa$ is $\aleph_1$ there is 
nothing to prove. 

Suppose that $\kappa>\aleph_1$ and \xitemof{chr-24} has been  
established for all infinite cardinals $<\kappa$. 

If $\kappa$ is a regular cardinal then 
\xitemof{chr-24} for $\kappa$ follows from \xitemof{chr-23}, the induction hypothesis and 
\Lemmaof{frp-2-0}. If $\kappa$ is a singular cardinal then \xitemof{chr-24} for 
$\kappa$ follows from \Thmof{Th-SC}, the induction hypothesis and \Lemmaof{frp-2-0}.
\qedof{\Thmof{main-thm-1}}

\label{literature}

\end{document}